\RequirePackage[l2tabu, orthodox]{nag}
\documentclass[a4paper,11pt]{article}
\usepackage{ifthen}
\newcommand{\OperationsResearchTemplate}{0}

\usepackage{authblk}
\usepackage{amsmath,amssymb,amsfonts}
\ifthenelse{\OperationsResearchTemplate = 0}{\usepackage{amsthm}}{}
\usepackage{bbm,booktabs}
\usepackage{color}
\usepackage{enumerate}
\usepackage{enumitem}
\ifthenelse{\OperationsResearchTemplate = 0}{\usepackage{fullpage}}{}
\usepackage[acronym]{glossaries}
\usepackage{graphicx}
\usepackage{multirow}
\usepackage{nicefrac}
\usepackage[percent]{overpic}
\usepackage{pgfplots}
\usepgfplotslibrary{groupplots}
\usepackage{setspace}
\usepackage{subfigure}
\usepackage{tikz}
\usetikzlibrary{arrows,automata,backgrounds,snakes}
\usepackage[colorinlistoftodos,textsize=scriptsize]{todonotes}
\usepackage{units}

\usepackage{epstopdf} % Needs to be loaded after loading graphicx.
\usepackage{hyperref}

\ifthenelse{\OperationsResearchTemplate = 0}{
\usepackage[T1]{fontenc}
\usepackage[bitstream-charter]{mathdesign}
}{}

\ifthenelse{\OperationsResearchTemplate = 0}{
}{
\OneAndAHalfSpacedXI % current default line spacing
\usepackage{endnotes}
\let\footnote=\endnote

\usepackage{natbib}
 \bibpunct[, ]{(}{)}{,}{a}{}{,}%
\TheoremsNumberedThrough     % Preferred (Theorem 1, Lemma 1, Theorem 2)
\ECRepeatTheorems

\EquationsNumberedThrough    % Default: (1), (2), ...
\MANUSCRIPTNO{todo} % When the article is logged in and DOI assigned to it, this manuscript number is no longer necessary
}

\definecolor{MyDarkBlue}{rgb}{0,0.08,0.50}
\definecolor{BrickRed}{rgb}{0.65,0.08,0}

\hypersetup{
colorlinks=true,       		% false: boxed links; true: colored links
    linkcolor=MyDarkBlue,   % color of internal links
    citecolor=BrickRed,     % color of links to bibliography
    filecolor=red,      	% color of file links
    urlcolor=cyan           % color of external links
}

\usepackage{dsfont} % For indicator function. \mathds type 1 font.

\renewcommand{\d}[1]{\ensuremath{\operatorname{d}\!{#1}}}
\newcommand{\e}[1]{ {\mathrm{e}}^{ #1 } }

\newcommand{\expectation}[1]{ \mathbb{E} [ #1 ] }

\newcommand{\indicator}[1]{ \mathds{1} [ #1 ] }

\newcommand{\process}[2]{ \{ #1 \}_{ #2 } }

\newcommand{\bigO}[1]{ O(#1) }
\newcommand{\bigObig}[1]{ O\Bigl(#1\Bigr) }

\newcommand{\probability}[1]{ \mathbb{P} [ #1 ] }

\newcommand{\naturalNumbersZero}{ \mathbb{N}_{0} }
\newcommand{\realNumbers}{ \mathbb{R} }
\newcommand{\positiveRealNumbers}{ [0,\infty) }

\newcommand{\criticalpoint}[1]{  #1^{\textnormal{opt}} }

\newcommand{\normalizationConstant}{Z}

\newcommand{\refFigure}[1]{{\textrm{Figure~\ref{#1}}}}

\newcommand{\refExample}[1]{{\textrm{Example~\ref{#1}}}}
\newcommand{\refEquation}[1]{{\textrm{\eqref{#1}}}}

\newcommand{\refProposition}[1]{{\textrm{Proposition~\ref{#1}}}}
\newcommand{\refLemma}[1]{{\textrm{Lemma~\ref{#1}}}}

\newcommand{\refSection}[1]{{\textrm{Section~\ref{#1}}}}
\newcommand{\refAppendixSection}[1]{\textrm{Appendix~\ref{#1}}}

\newcommand{\QuodEratDemonstrandum}{\hfill \ensuremath{\Box}}

\ifthenelse{\OperationsResearchTemplate = 0}{

\newtheorem{proposition}{Proposition}
\newtheorem{lemma}{Lemma}
\newtheorem{corollary}{Corollary}

\theoremstyle{definition}
\newtheorem{example}{Example}
\theoremstyle{remark}

}{}

\begin{document}

\ifthenelse{\OperationsResearchTemplate = 1}
{\newcommand\myTitle[1]{\RUNTITLE{Optimal Admission Control for Many-Server Systems with QED-Driven Revenues}\TITLE{Optimal Admission Control for Many-Server Systems \\ with QED-Driven Revenues}}}
{\newcommand\myTitle[1]{\title{Optimal Admission Control for Many-Server Systems \\ with QED-Driven Revenues}}}

\ifthenelse{\OperationsResearchTemplate = 1}
{\newcommand\myAbstract[1]{\ABSTRACT{#1}}}
{\newcommand\myAbstract[1]{\begin{abstract}#1\end{abstract}}}

\ifthenelse{\OperationsResearchTemplate = 1}
{\newcommand\myProof[1]{\proof{Proof}{#1 \QuodEratDemonstrandum}}}
{\newcommand\myProof[1]{\begin{proof}#1\end{proof}}}

\ifthenelse{\OperationsResearchTemplate = 1}
{\newcommand\myAcknowledgment[1]{\ACKNOWLEDGMENT{#1}}}
{\newcommand\myAcknowledgment[1]{\section*{Acknowledgment}#1}}

\ifthenelse{\OperationsResearchTemplate = 1}
{\newcommand\myAppendices[1]{\ECSwitch \ECHead{E-companion} #1}}
{\newcommand\myAppendices[1]{\appendix #1}}

\makeglossaries
\newacronym{QED}{QED}{Quality-and-Efficiency-Driven}
\newacronym{OU}{OU}{Ornstein--Uhlenbeck}
\newacronym{EM}{EM}{Euler--Maclaurin}

\myTitle{}

%%%%%%%%%%%%%%%%%%%%%%%%%%%
%%% Author information. %%%
%%%%%%%%%%%%%%%%%%%%%%%%%%%
\ifthenelse{\OperationsResearchTemplate = 1}
{
\RUNAUTHOR{Jaron Sanders et al.}
\ARTICLEAUTHORS{%
\AUTHOR{Jaron Sanders}
\AFF{\EMAIL{jaron.sanders@tue.nl}} %, \URL{}}
\AUTHOR{S.C.\ Borst}
\AFF{\EMAIL{s.c.borst@tue.nl}}
\AUTHOR{A.J.E.M.\ Janssen}
\AFF{\EMAIL{a.j.e.m.janssen@tue.nl}}
\AUTHOR{J.S.H.\ van Leeuwaarden}
\AFF{\EMAIL{j.s.h.v.leeuwaarden@tue.nl} \\ \vspace{1em} Dept.\ of Mathematics \& Computer Science, Eindhoven University of Technology, The Netherlands}
}%
}{
\date{\today}
\author{Jaron Sanders}
\author{S.C.\ Borst}
\author{A.J.E.M.\ Janssen}
\author{J.S.H.\ van Leeuwaarden\footnote{E-mail addresses: \texttt{jaron.sanders@tue.nl}, \texttt{s.c.borst@tue.nl}, \texttt{a.j.e.m.janssen@tue.nl}, and \texttt{j.s.h.v.leeuwaarden@tue.nl}}}
\affil{Department of Mathematics \& Computer Science, Eindhoven University of Technology, P.O.\ Box 513, 5600 MB Eindhoven, The Netherlands}
}

%%%%%%%%%%%%%%
%%% Title. %%%
%%%%%%%%%%%%%%
\ifthenelse{\OperationsResearchTemplate = 0}{\maketitle}{}

\myAbstract{
We consider Markovian many-server systems with admission control operating in a \gls{QED} regime, where the relative utilization approaches unity while the number of servers grows large, providing natural Economies-of-Scale. In order to determine the optimal admission control policy, we adopt a revenue maximization framework, and suppose that the revenue rate attains a maximum when no customers are waiting and no servers are idling. When the revenue function scales properly with the system size, we show that a nondegenerate optimization problem arises in the limit. Detailed analysis demonstrates that the revenue is maximized by nontrivial policies that bar customers from entering when the queue length exceeds a certain threshold of the order of the typical square-root level variation in the system occupancy. We identify a fundamental equation characterizing the optimal threshold, which we extensively leverage to provide broadly applicable upper/lower bounds for the optimal threshold, establish its monotonicity, and examine its asymptotic behavior, all for general revenue structures. For linear and exponential revenue structures, we present explicit expressions for the optimal threshold.
\ifthenelse{\OperationsResearchTemplate = 0}{
}{}
}

\ifthenelse{\OperationsResearchTemplate = 1}{
}{}

%%%%%%%%%%%%%%
%%% Title. %%%
%%%%%%%%%%%%%%
\ifthenelse{\OperationsResearchTemplate = 1}{\maketitle}{}

%%%%%%%%%%%%%%%%%%%%%%%%%%%%%%%
%%% Introductie versie Sem. %%%
%%%%%%%%%%%%%%%%%%%%%%%%%%%%%%%
\section{Introduction}
\label{sec:Introduction}

Large-scale systems that operate in the Quality-and-Efficiency Driven (QED)
regime dwarf the usual trade-off between high system utilization and short
waiting times.
In order to achieve these dual goals, the system is scaled so as to approach
full utilization, while the number of servers grows simultaneously large,
rendering crucial Economies-of-Scale.
Specifically, for a Markovian many-server system with Poisson arrival
rate $\lambda$, exponential unit-mean service times and $s$ servers, the load
$\rho = \lambda / s$ is driven to unity in the \gls{QED} regime in accordance with
\begin{equation}
\label{1.1}
(1 - \rho) \sqrt{s} \to \gamma, \quad s \to \infty,
\end{equation}
for some fixed parameter $\gamma \in \mathbb{R}_+$.
As $s$ grows large, the stationary probability of delay then tends to a limit, say $g(\gamma)$, which may take any value in $(0, 1)$, depending on the
parameter $\gamma$.
Since the \textit{conditional} queue length distribution is geometric with mean
$\rho / (1 - \rho) \approx \sqrt{s} / \gamma$, it follows that the stationary
mean number of waiting customers scales as $g(\gamma) \sqrt{s} / \gamma$.
Little's law then in turn implies that the mean stationary waiting time
of a customer falls off as $g(\gamma) / ( \gamma \sqrt{s} )$.

The \gls{QED} scaling behavior also manifests itself in process level limits,
where the evolution of the system occupancy, properly centered around $s$
and normalized by $\sqrt{s}$, converges to a diffusion process as $s\to\infty$,
which again is fully characterized by the single parameter $\gamma$.
This reflects that the system state typically hovers around the full-occupancy
level $s$, with natural fluctuations of the order $\sqrt{s}$.

The \gls{QED} scaling laws provide a powerful framework for system dimensioning,
i.e., matching the service capacity and traffic demand so as to achieve
a certain target performance level or optimize a certain cost metric.
Suppose, for instance, that the objective is to find the number of servers $s$
for a given arrival rate $\lambda$ (or equivalently, determine what arrival
rate $\lambda$ can be supported with a given number $s$ of servers)
such that a target delay probability $\epsilon \in (0, 1)$ is attained.
The above-mentioned convergence results for the delay probability then
provide the natural guideline to match the service capacity and traffic
volume in accordance with $\lambda = s - \gamma_\epsilon \sqrt{s}$, where
the value of $\gamma_\epsilon$ is such that $g(\gamma_\epsilon) = \epsilon$.

As an alternative objective, imagine we aim to strike a balance between
the expenses incurred for staffing servers and the dissatisfaction
experienced by waiting customers.
Specifically, suppose a (salary) cost $c$ is associated with each server
per unit of time and a (possibly fictitious) holding charge $h$ is imposed
for every waiting customer per unit of time.
Writing $\lambda = s - \gamma \sqrt{s}$ in accordance with \refEquation{1.1},
and recalling that the mean number of waiting customers scales as
$g(\gamma) \sqrt{s} / \gamma$, we find that the total operating cost per
time unit scales as
\begin{equation}
c s + h \frac{g(\gamma) \sqrt{s}}{\gamma} =
\lambda c + c \gamma \sqrt{s} + h \frac{g(\gamma) \sqrt{s}}{\gamma} =
\lambda c + \left( c \gamma + h \frac{g(\gamma)}{\gamma}\right) \sqrt{s}.
\end{equation}
This then suggests to set the number of servers in accordance with
$s = \lambda + \gamma_{c,h} \sqrt{s}$, where
$\gamma_{c,h} = \arg \min_{\gamma > 0} (c \gamma + h g(\gamma) / \gamma )$
in order to minimize the total operating cost per time unit.
Exploiting the powerful \gls{QED} limit theorems, such convenient capacity sizing rules
can in fact be shown to achieve optimality in some suitable asymptotic sense.

As illustrated by the above two paragraphs, the \gls{QED} scaling laws can be
leveraged for the purpose of dimensioning, with the objective to balance
the service capacity and traffic demand so as to achieve a certain target
performance standard or optimize a certain cost criterion.
A critical assumption, however, is that all customers are admitted into the
system and eventually served, which may in fact not be consistent with the
relevant objective functions in the dimensioning, let alone be optimal in
any sense.

Motivated by the latter observation, we focus in the present paper on the
optimal admission control problem for a given performance or cost criterion.
Admission control acts on an operational time scale, with decisions
occurring continuously whenever customers arrive, as opposed to capacity
planning decisions which tend to involve longer time scales.
Indeed, we assume that the service capacity and traffic volume are given,
and balanced in accordance with \refEquation{1.1}, but do allow for the value
of $\gamma$ to be negative, since admission control provides a mechanism
to deal with overload conditions.
While a negative value of $\gamma$ may not be a plausible outcome of
a deliberate optimization process, in practice an overload of that order
might well result from typical forecast errors.

We formulate the admission control problem in a revenue maximization
framework, and suppose that revenue is generated at rate $r_s(k)$ when the
system occupancy is $k$.
As noted above, both from a customer satisfaction perspective and a system
efficiency viewpoint, the ideal operating condition for the system is
around the full occupancy level $s$, where no customers are waiting and no
servers are idling.
Hence we assume that the function $r_s(k)$ is unimodal, increasing
in $k$ for $k \leq s$ and decreasing in $k$ for $k \geq s$, thus attaining
its maximum at $k = s$.

We consider probabilistic control policies, which admit arriving customers
with probability $p_s(k-s)$ when the system occupancy is $k$,
independent of any prior admission decisions.
It is obviously advantageous to admit customers as long as free servers are
available, since it will not lead to any wait and drive the system closer
to the ideal operating point $s$, boosting the instantaneous revenue rate.
Thus we stipulate that $p_s(k-s) = 1$ for all $k < s$.

For $k \geq s$, it is far less evident whether to admit customers or not.
Admitting a customer will then result in a wait and move the system away
from the ideal operating point, reducing the instantaneous revenue rate.
On the other hand, admitting a customer may prevent the system occupancy
from falling below the ideal operating point in the future.
The potential long-term gain may outweigh the adverse near-term effect,
so there may be a net benefit, but the incentive weakens as the queue grows.
The fundamental challenge in the design of admission control policies is
to find exactly the point where the marginal utility reaches zero,
so as to strike the optimal balance between the conflicting near-term
and longer-term considerations.

Since the service capacity and traffic volume are governed by \refEquation{1.1},
the \gls{QED} scaling laws imply that, at least for $\gamma > 0$ and without any
admission control, the system occupancy varies around the ideal operating
point $s$, with typical deviations of the order $\sqrt{s}$.
It is therefore natural to suppose that the revenue rates and admission
probabilities scale in a consistent manner, and in the limit behave as functions
of the properly centered and normalized state variable $(k - s) / \sqrt{s}$. Specifically, we assume that the revenue rates satisfy the scaling condition
\begin{equation}
\label{formr}
\frac{r_s(k) - n_s}{q_s} \to r\Bigl(\frac{k-s}{\sqrt{s}}\Bigr),
\quad s \to \infty,
\end{equation}
with $n_s$ a nominal revenue rate attained at the ideal operating point,
$q_s$ a scaling coefficient, and $r$ a unimodal function,
which represents the scaled reduction in revenue rate associated with
deviations from the optimal operating point $s$.
For example, with $[x]^+ = \max\{ 0, x \}$, any revenue structure of the form
\begin{equation}
r_s(k) = n_s - \alpha^- ([s-k]^+)^{\beta^-} - \alpha^+ ([k-s]^+)^{\beta^+}
\end{equation}
satisfies \refEquation{formr} when $q_s = s^{\max\{\beta^-, \beta^+\} / 2}$, in which case
\begin{equation}
r(x) = - \alpha^- ([- x]^+)^{\beta^-} \indicator{ \beta^- \geq \beta^+ }
- \alpha^+ ([x]^+)^{\beta^+} \indicator{ \beta^- \leq \beta^+ }.
\end{equation}
Note that these revenue structures impose
polynomial penalties on deviations from the ideal operating point. Similar to \refEquation{formr}, we assume that the admission probabilities satisfy a scaling condition, namely
\begin{equation}
\label{formf}
p_s(0) \cdots p_s(k-s) = f\Bigl(\frac{k-s}{\sqrt{s}}\Bigr),
\quad k \geq s,
\end{equation}
with $f$ a non-increasing function and $f(0) = 1$.
In particular, we allow for $f(x) = \indicator{0 \leq x < \eta}$, which 
corresponds to an admission threshold control $p_s(k-s) = \indicator{ k -s \leq \lfloor \eta \sqrt{s} \rfloor }$.

In \refSection{sec:Revenue_framework} we discuss the fact that the optimal admission policy is indeed such a threshold
control, with the value of $\eta$ asymptotically being determined by the
function $r$, which we later prove in \refSection{sec:Optimality_of_threshold_policies}.
The optimality of a threshold policy may not come as a surprise,
and can in fact be established in the pre-limit ($s < \infty$).
However, the pre-limit optimality proof only yields the structural property,
and does not furnish any characterization of how the optimal threshold depends
on the system characteristics or provide any computational procedure for
actually obtaining the optimal value.
In contrast, our asymptotic framework (as $s \to \infty$) produces a specific equation
characterizing the optimal threshold value, which does offer explicit
insight in the dependence on the key system parameters and can serve as
a basis for an efficient numerical computation or even a closed-form
expression in certain cases.
This is particularly valuable for large-scale systems where a brute-force
enumerative search procedure may prove prohibitive.

Let us finally discuss the precise form of the revenue rates that serve as the objective function that needs to be maximized by the optimal threshold. We will mostly focus on the average \emph{system-governed} revenue rate defined as
\begin{equation}
R_s(\{p_s(k)\}_{k\geq0}) 
= \sum_{k = 0}^{\infty} r_s(k) \pi_s(k).
\label{eqn:System_governed_revenue_rate}
\end{equation}
From the system's perspective, this means that the revenue is simply governed by the state-dependent revenue rate $r_s(k)$ weighed according to the stationary distribution, with $\pi_s(k)$ denoting the stationary probability of state $k$.

An alternative would be to consider the \emph{customer reward} rate 
\begin{equation}
\hat{R}_s(\{p_s(k)\}_{k\geq0}) 
= \lambda \sum_{k = 0}^{\infty} \hat{r}_s(k) p_s(k-s) \pi_s(k).
\label{eqn:Customer_oriented_reward_rate}
\end{equation}
Here, $\hat{r}_s(k)$ can be interpreted as the state-dependent reward when admitting a customer in state $k$, and since this happens with probability $p_s(k)$ at intensity $\lambda$, we obtain \refEquation{eqn:Customer_oriented_reward_rate}. While this paper primarily focuses on \refEquation{eqn:System_governed_revenue_rate}, we show in \refSection{sec:Revenue_framework__Customer_state_dependent_rewards} that there is an intimate connection with \refEquation{eqn:Customer_oriented_reward_rate}; a system-governed reward structure $\{ r_s(k) \}_{k \in \naturalNumbersZero}$ can be translated into a customer reward structure $\{ \hat{r}_s(k) \}_{k \in \naturalNumbersZero}$, and vice versa.

%%%%%%%%%%%%%%%%%%%%%%%%%%%%%%
%%% Literature survey Sem. %%%
%%%%%%%%%%%%%%%%%%%%%%%%%%%%%%
\subsection{Contributions and related literature}

% Optimality of thresholds.
A diverse range of control problems have been considered in the queueing
literature, and we refer the reader to
\cite{lippman_applying_1975,stidham_optimal_1985,kushner_numerical_2001,meyn_control_2008,cil_effects_2009}
for background.
Threshold control policies are found to be optimal in a variety of contexts such as \cite{de_waal_overload_1990,chen_state_2001,bekker_optimal_2006},
and many (implicit) characterizations of optimal threshold values have been
obtained in \cite{naor_regulation_1969,yildirim_admission_2010}, and \cite{borgs_optimal_2014}.
For (single-server) queues in a conventional heavy-traffic regime,
optimality of threshold control policies has been established by studying
limiting diffusion control problems in \cite{ghosh_optimal_2007,ward_asymptotically_2008}, and \cite{ghosh_optimal_2010}.

%Brownian motion in control.
%Our focus is on identifying optimal admission controls for many-server queueing systems that operate in the balanced heavy-traffic regime known as the \gls{QED} or Halfin-Whitt regime \cite{halfin_heavy-traffic_1981}. This follows recent advances on appropriate incorporation and scaling of admission controls in the \gls{QED} regime \cite{janssen_scaled_2013}, with which relevant effect is achieved on the limiting diffusion processes. Note that Brownian control in general has many applications \cite{harrison_brownian_2013}.

% Optimal admission control in the \gls{QED} regime
The analysis of control problems in the \gls{QED} regime has mostly focused
on routing and scheduling, see \cite{atar_brownian_2004,atar_scheduling_2005,atar_diffusion_2005,atar_queueing_2006}, and \cite{gurvich_queue-and-idleness-ratio_2009}. Threshold policies in the context of many-server systems in the \gls{QED} regime have been considered in \cite{armony_customer_2004,massey_asymptotically_2005,whitt_heavy-traffic_2005}, and \cite{whitt_diffusion_2004}.
General admission control, however, has only received limited attention in the \gls{QED} regime,
see for instance \cite{kocaga_admission_2010,weerasinghe_abandonment_2013}.
These studies specifically account for abandonments, which create
a trade-off between the rejection of a new arrival and the risk of that
arrival later abandoning without receiving service, with the associated
costly increase of server idleness.

% Appropriate scaling.
In the present paper we address the optimal admission control problem
from a revenue maximization perspective.
We specifically build on the recent work in \cite{janssen_scaled_2013}
to show that a nondegenerate optimization problem arises in the limit
when the revenue function scales properly with the system size.
Analysis of the latter problem shows that nontrivial threshold control
policies are optimal in the \gls{QED} regime for a broad class of revenue
functions that peak around the ideal operating point.

%The difficulty in \gls{QED} admission control problems is that in order to obtain non-trivial optimal admission control policies, \emph{all} system parameters need to be properly scaled. In \refSection{sec:Revenue_framework}, we build on the recent work in \cite{janssen_scaled_2013} and design a framework that scales the revenue structure, control policy, and load, appropriately, so that a nondegenerate limiting control problem results. If for example the load is not suitably scaled, the asymptotic optimal control policy tends to one of the extreme policies that never admit or never reject a customer. 
%This emphasizes the importance of designing a framework that incorporates appropriate scaling. As far as we have been able to assert, this has not yet been discussed in literature. 

% Implicit characterizations.
In \refSection{sec:Revenue_framework} we present a fundamental equation
which implicitly determines the asymptotically optimal threshold.
The subsequent analysis of this equation in
\refSection{sec:Analysis_of_optimal_thresholds} yields valuable insight into
the dependence of the optimal threshold on the revenue structure,
and provides a basis for an efficient numerical scheme.
Closed-form expressions for the optimal threshold can only be derived when considering
specific revenue structures.

We will, for example, show that for \emph{linearly decreasing} revenue rates,
the optimal threshold can be (explicitly) expressed in terms of the Lambert W
function \cite{corless_lambert_1996}. We note that a linearly decreasing revenue structure has also been considered in \cite{borgs_optimal_2014} for determining the optimal threshold $\criticalpoint{k}$ in an $\mathrm{M}/\mathrm{M}/s/k$ system, and there also, $\criticalpoint{k}$ is expressed in terms of the Lambert W function. Besides assuming a static revenue and finite threshold $k$, a crucial difference between \cite{borgs_optimal_2014} and this paper is that our revenue structure scales as in \refEquation{formr}, so that the threshold $k$ is suitable for the \gls{QED} regime. Our work thus extends \cite{borgs_optimal_2014}, both in terms of scalable and more general revenue structures.

% Mathematical techniques and their use in queueing \gls{QED} literature.
In terms of mathematical techniques, we use \gls{EM} summation
(\cite{olver_NIST_2010}) to analyze the asymptotic behavior of \refEquation{eqn:System_governed_revenue_rate} as $s \to \infty$.
This approach was used recently for many-server systems with admission
control in the \gls{QED} regime \cite{janssen_scaled_2013},
and is now extended by incorporating suitably scaled revenue structures
in \refSection{sec:Revenue_framework}.
These ingredients then pave the way to determine the optimal admission
control policy in the \gls{QED} regime in
\refSection{sec:Analysis_of_optimal_thresholds}.
In \refSection{sec:Optimality_of_threshold_policies},
we use Hilbert-space theory from analysis,
and techniques from variational calculus,
to prove the existence of optimal control policies, and to establish that
control policies with an admission threshold which scales with the natural
$\sqrt{s}$ order of variation are optimal in the \gls{QED} regime.

\section{Revenue maximization framework}
\label{sec:Revenue_framework}

\glsresetall

We now develop an asymptotic framework for determining an optimal
admission control policy for a given performance or cost criterion.
In \refSection{sec:Revenue_framework__MMs_type_queue_with_admission_control}
we describe the basic model for the system dynamics, which is an extension
of the classical $\mathrm{M}/\mathrm{M}/s$ system.
Specifically, the model incorporates admission control and is augmented
with a revenue structure, which describes the revenue rate as a function of
the system occupancy.
Adopting this flexible apparatus, the problem of finding an optimal admission
control policy is formulated in terms of a revenue maximization objective.

\subsection{Markovian many-server systems with admission control}
\label{sec:Revenue_framework__MMs_type_queue_with_admission_control}

\begin{figure}[!hbtp]
\begin{center}
\small
\begin{tikzpicture}[ ->, >=stealth, shorten >=1pt, auto, node distance=2.5cm, semithick, sloped]
\tikzstyle{every state} = [ fill=none, thick, scale=0.8, minimum height=1.25cm]
\node[state] (s0) {$0$};
\node[state] (s1) [right of=s0] {$1$};
    \path (s0) edge  [bend right=10] node[below] {$\lambda$} (s1);
    \path (s1) edge  [bend right=10] node[above] {$1$} (s0);
\node[state] (s2) [right of=s1, draw=none] {$\cdots$};
    \path (s1) edge  [bend right=10] node[below] {$\lambda$} (s2);
    \path (s2) edge  [bend right=10] node[above] {$2$} (s1);
\node[state] (s4) [right of=s2] {$s$};    
    \path (s2) edge  [bend right=10] node[below] {$\lambda$} (s4);
    \path (s4) edge  [bend right=10] node[above] {$s$} (s2);
\node[state] (s5) [right of=s4] {$s+1$};    
    \path (s4) edge  [bend right=10] node[below] {$p_s(0) \lambda$} (s5);
    \path (s5) edge  [bend right=10] node[above] {$s$} (s4);        
\node[state] (s6) [right of=s5, draw=none] {$\cdots$};    
    \path (s5) edge  [bend right=10] node[below] {$p_s(1) \lambda$} (s6);
    \path (s6) edge  [bend right=10] node[above] {$s$} (s5);    
\end{tikzpicture}
\end{center}
\vspace{-1em}
\caption{\textrm{Transition diagram of the process $\process{ Q_s(t) }{ t \geq 0 }$.}}
\label{fig:Figure__Markov_chain_of_admission_controlled_process}
\end{figure}
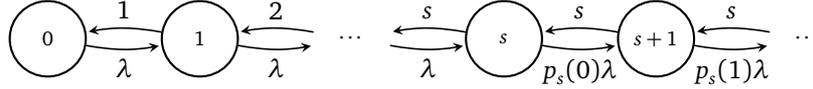

Consider a system with $s$ parallel servers where customers arrive
according to a Poisson process with rate $\lambda$.
Customers require exponentially distributed service times with unit mean.
A customer that finds upon arrival $k$ customers in the system is taken
into service immediately if $k < s$, or may join a queue of waiting
customers if $k \geq s$.
If all servers are occupied, a newly arriving customer is admitted into the
system with probability $p_s(k - s)$, and denied access otherwise.
We refer to the probabilities $\{p_s(k)\}_{k \geq 0}$ as the
\emph{admission control policy}.
If we denote the number of customers in the system at time $t$ by $Q_s(t)$,
and make the usual independence assumptions, then $\process{Q_s(t)}{t \geq 0}$
constitutes a Markov process (see \refFigure{fig:Figure__Markov_chain_of_admission_controlled_process} for its transition diagram).
The stationary distribution
$\pi_s(k) = \lim_{t \rightarrow \infty} \probability{Q_s(t) = k}$ is given by
\begin{equation}
\pi_s(k)
=
\begin{cases}
\pi_s(0) \frac{(s\rho)^k}{k!}, & k = 1, 2, \ldots, s, \\
\pi_s(0) \frac{ s^s \rho^k }{ s! } \prod_{i=0}^{k-s-1} p_s(i),  & k = s+1, s+2, \ldots,
\end{cases}
\label{eqn:Equilibrium_distribution}
\end{equation}
with
\begin{equation}
\rho
= \frac{\lambda}{s},
\quad
\pi_s(0) = \Bigl( \sum_{k=0}^s\,\frac{(s\rho)^k}{k!}+\frac{(s\rho)^s}{s!}F_s(\rho) \Bigr)^{-1} \label{eqn:Normalization_constant}
\end{equation}
and
\begin{equation}
F_s(\rho)=\sum_{n=0}^{\infty}\,p_s(0) \cdots p_s(n)\,\rho^{n+1}. \label{eqn:Fs}
\end{equation}
 From \refEquation{eqn:Equilibrium_distribution}--\refEquation{eqn:Fs},
we see that the stationary distribution exists if and only if the relative
load $\rho$ and the admission control policy
$\{p_s(k)\}_{k \in \naturalNumbersZero}$ are such that $F_s(\rho) < \infty$
\cite{janssen_scaled_2013}, which always holds in case $\rho < 1$.

With $k$ customers in the system, we assume that the system generates revenue at
rate $r_s(k) \in \realNumbers$. We call
$\{r_s(k)\}_{k \geq 0}$ the \emph{revenue structure}.
Our objective is to find an admission control policy in terms of the
probabilities $\{p_s(k)\}_{k \geq 0}$ that maximizes the average stationary
revenue rate, i.e.\

\begin{equation}
\begin{array}{ccccc}
\textrm{to maximize} & R_s(\{p_s(k)\}_{k\geq0}) & \textrm{over} & \{ p_s(k) \}_{ k \geq 0}, \\
\textrm{subject to} & 0 \leq p_s(k) \leq 1, \quad k \in \naturalNumbersZero, & \textrm{and} & F_s(\rho) < \infty. \\
\end{array}
\label{eqn:Objective_optimization_problem_prelimit}
\end{equation}

\subsection{QED-driven asymptotic optimization framework}
\label{sec:Revenue_framework__Asymptotic_control_problem}

We now construct an asymptotic optimization framework where
the limit laws of the \gls{QED} regime can be leveraged by imposing
suitable assumptions on the admission control policy and revenue structure.
In order for the system to operate in the \gls{QED} regime, we couple the
arrival rate to the number of servers as
\begin{equation}
\lambda = s - \gamma \sqrt{s},
\quad \gamma \in \realNumbers.
\label{eqn:QED_scaling}
\end{equation}
For the admission control policy we assume the form in \refEquation{formf},
with $f$ either a nonincreasing, bounded, and twice differentiable
continuous function, or a step function, which we will refer to as the \emph{asymptotic
admission control profile}.
We also assume the revenue structure has the scaling property \refEquation{formr},
with $r$ a piecewise bounded, twice differentiable continuous
function with bounded derivatives. We will refer to $r$ as the
\emph{asymptotic revenue profile}.
These assumptions allow us to establish
\refProposition{prop:QED_limit_of_revenue__Continuous_admission_controls}
by considering the stationary average revenue rate
$R_s(\{p_s(k)\}_{k\geq0})$ as a Riemann sum and using
\gls{EM} summation to identify its limiting integral expression,
the proof of which can be found in
\refAppendixSection{sec:Limiting_behavior_of_long_term_QED_revenue}. Let
$
\phi(x) = \exp{( - \frac{1}{2} x^2 )} / \sqrt{2\pi}
$ and
$
\Phi(x) = \int_{-\infty}^x \phi(u) \d{u}
$
denote the probability density function and cumulative distribution function of the standard normal distribution, respectively.

\begin{proposition}
\label{prop:QED_limit_of_revenue__Continuous_admission_controls}
If $r^{(i)}$ is continuous and bounded for $i = 0$, $1$, $2$, and either
$\mathrm{(i)}$ $f$ is smooth, and $(f(x) \exp{(- \gamma x)})^{(i)}$
is exponentially small as $x \rightarrow \infty$ for $i = 0$, $1$, $2$,
or $\mathrm{(ii)}$ $f(x) = \indicator{0 \leq x < \eta}$ with a fixed,
finite $\eta > 0$, then
\begin{equation}
\lim_{s \rightarrow \infty} \frac{ R_s(\{p_s(k)\}_{k\geq0}) - n_s }{q_s} = R(f),
\end{equation}
where in case $\mathrm{(i)}$
\begin{equation}
R(f)
= \frac{ \int_{-\infty}^0 r(x) \e{ - \frac{1}{2} x^2 - \gamma x } \d{x} + \int_0^\infty r(x) f(x) \e{-\gamma x} \d{x} }{ \frac{\Phi(\gamma)}{\phi(\gamma)} + \int_0^\infty f(x) \e{-\gamma x} \d{x} },
\label{eqn:QED_limit_of_revenue_rate}
\end{equation}
and in case $\mathrm{(ii)}$
\begin{equation}
R(\indicator{0 \leq x < \eta})
= \frac{ \int_{-\infty}^0 r(x) \e{- \frac{1}{2} x^2 - \gamma x } \d{x} + \int_0^\eta r(x) \e{-\gamma x} \d{x} }{ \frac{\Phi(\gamma)}{\phi(\gamma)} + \frac{1 - \e{-\gamma \eta}}{\gamma} },
\quad \eta \geq 0.
\label{eqn:Long_term_revenue__Threshold}
\end{equation}
\end{proposition}

Because of the importance of the threshold policy, we will henceforth use the short-hand notations $R_{\mathrm{T},s}(\tau) = R_s( \{ \indicator{k \leq s + \tau} \}_{k\geq0} )$ and $R_{\mathrm{T}}(\eta) 
= R(\indicator{0 \leq x < \eta})$ to indicate threshold policies.

\begin{example}[Exponential revenue]
\label{expo}
Consider a revenue structure $r_s(k) = \exp{(b (k - s) / \sqrt{s})}$ for $k < s$
and $r_s(k) = \exp{( - d (k - s) / \sqrt{s} )}$ for $k \geq s$, and with $b, d > 0$. Taking $n_s = 0$ and $q_s = 1$, the asymptotic revenue profile is
$r(x) = \exp{(bx)}$ for $x < 0$ and $r(x) = \exp{(-dx)}$ for $x \geq 0$, so that according to \refProposition{prop:QED_limit_of_revenue__Continuous_admission_controls} for threshold policies,
\begin{equation}
\lim_{s \to \infty} R_{\mathrm{T},s}( \lfloor \eta \sqrt{s} \rfloor )
= R_{\mathrm{T}}(\eta) 
= \frac{ \frac{\Phi(\gamma-b)}{\phi(\gamma-b)} + \frac{ 1 - \e{-(d+\gamma)\eta} }{d+\gamma} }{ \frac{\Phi(\gamma)}{\phi(\gamma)} + \frac{1 - \e{-\gamma \eta}}{\gamma} }.
\end{equation}
\end{example}

\begin{figure}[!hbtp]
\begin{center}
%\subfigure{
%\includegraphics[width=0.8\columnwidth]{Figure__Revenue_structure_as_s_increases}
\small
\input{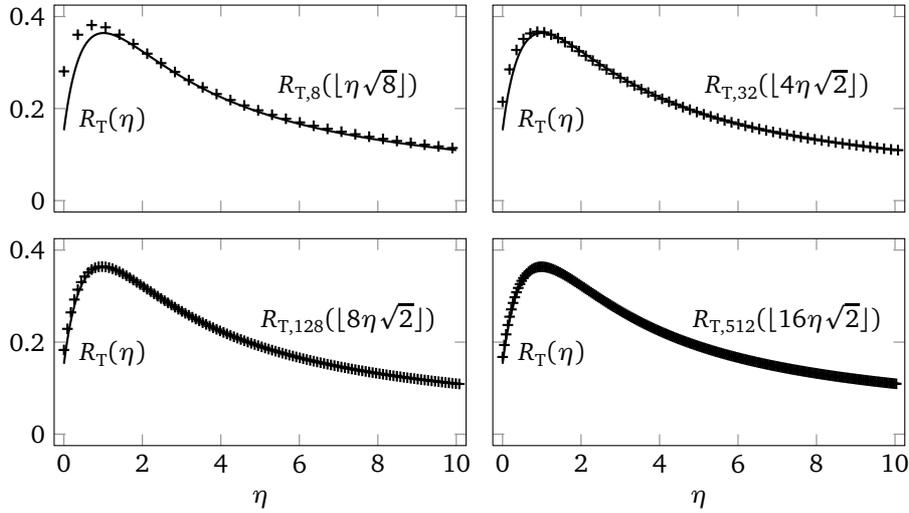}
%}
\end{center}
\vspace{-1em}
\caption{\textrm{$R_{\mathrm{T},s}( \lfloor \eta \sqrt{s} \rfloor )$ and $R_{\mathrm{T}}(\eta)$ for $s = 8$, $32$, $128$, $256$ servers.}}
\label{fig:Figure__Revenue_structure_as_s_increases}
\end{figure}

\refFigure{fig:Figure__Revenue_structure_as_s_increases} plots $R_{\mathrm{T},s}( \lfloor \eta \sqrt{s} \rfloor )$
for a finite system with $s = 8$, $32$, $128$, $256$ servers, respectively,
together with its limit $R_{\mathrm{T}}(\eta)$. Here, we set $b = 5$, $d = 1$, and $\gamma = 0.01$.
Note that the approximation
$ R_{\mathrm{T},s}( \lfloor \eta \sqrt{s} \rfloor ) \approx R_{\mathrm{T}}(\eta)$
is remarkably accurate, even for relatively small systems, an observation
which in fact seems to hold for most revenue structures and parameter choices. For this particular revenue structure, we see that the average revenue rate
peaks around $\criticalpoint{\eta} \approx 1.0$. In \refExample{eg:Exponential_revenue_revisited}, we confirm this observation by determining $\criticalpoint{\eta}$ numerically.

An alternative way of establishing that the limit of $R_s(\{p_s(k)\}_{k \geq 0})$ is $R(f)$,
is by exploiting the stochastic-process limit for $\process{Q_s(t)}{t \geq 0}$.
It was shown in \cite{janssen_scaled_2013} that under condition~(i)
in \refProposition{prop:QED_limit_of_revenue__Continuous_admission_controls}, together with \refEquation{formf} and \refEquation{eqn:QED_scaling}, the normalized process
$\process{\hat{Q}_s(t)}{t \geq 0}$ with $\hat{Q}_s(t) = (Q_s(t) - s) / \sqrt{s}$
converges weakly to a stochastic-process limit $\process{D(t)}{t \geq 0}$
with stationary density 
\begin{equation}
w(x) 
= 
\begin{cases}
\normalizationConstant^{-1} \e{ - \frac{1}{2} x^2 - \gamma x }, & x < 0, \\
\normalizationConstant^{-1} f(x) \e{- \gamma x}, & x \geq 0, \\
\end{cases}
\end{equation}
where $\normalizationConstant = \Phi(\gamma) / \phi(\gamma) + \int_0^\infty f(x) \exp{(- \gamma x)} \d{x}$. When additionally assuming \refEquation{formr}, the limiting system revenue can be written as 
\begin{equation}
R(f) = \int_{- \infty}^\infty r(x) w(x) \d{x},
\end{equation}
the stationary revenue rate generated by the stochastic-process limit. So an alternative method to prove \refProposition{prop:QED_limit_of_revenue__Continuous_admission_controls} would be to 
first formally establish weak convergence at the process level, then prove that limits with respect to space and time can be interchanged, and finally use the
stationary behavior of the stochastic-process limit. This is a common approach in the \gls{QED} literature (\cite{halfin_heavy-traffic_1981,garnett_designing_2002}). Instead, we construct a direct, purely analytic proof, that additionally gives insight into 
%the rate at which the limit $R(f)$ is reached, or in other words,
the error that is made when approximating $R_s(\{p_s(k)\}_{k\geq0})$ by $R(f)$ for finite $s$.
These error estimates are available in \refAppendixSection{sec:Limiting_behavior_of_long_term_QED_revenue} for future reference.

With \refProposition{prop:QED_limit_of_revenue__Continuous_admission_controls} at hand, we are naturally led to consider the asymptotic optimization problem, namely,
\begin{equation}
\begin{array}{ccccc}
\textrm{to maximize} & R(f) & \textrm{over} & f, \\
\textrm{subject to} &  0 \leq f(x) \leq 1, \quad x \in \positiveRealNumbers, & \textrm{and} & \int_{0}^\infty f(x) \e{-\gamma x} \d{x} < \infty. \\
\end{array}
\label{eqn:Objective_optimization_problem_postlimit}
\end{equation}
The condition $\int_0^\infty f(x) \e{- \gamma x} \d{x} < \infty$ is
the limiting form of the stability condition $F_s(\rho) < \infty$, see \cite{janssen_scaled_2013}. Also note that we do not restrict $f$ to be monotone. We prove for the optimization problem in \refEquation{eqn:Objective_optimization_problem_postlimit} the following in \refSection{sec:Optimality_of_threshold_policies}.

\begin{proposition}
\label{prop:Threshold_control_is_optimal_form_of_control}
If $r$ is nonincreasing for $x \geq 0$, then there exist optimal
asymptotic admission controls that solve \refEquation{eqn:Objective_optimization_problem_postlimit}.
Moreover, the optimal asymptotic admission control profiles have
a threshold structure of the form
\begin{equation}
f(x) = \indicator{ 0 \leq x < \criticalpoint{\eta} },
\end{equation} 
where $\criticalpoint{\eta}$ is any solution of
\begin{equation}
r(\eta) = R_{\mathrm{T}}(\eta) \label{eqn:Threshold_equation}
\end{equation}
if $r(0) > R_{\mathrm{T}}(0)$, and $\criticalpoint{\eta} = 0$
if $r(0) \leq R_{\mathrm{T}}(0)$.
If $r$ is strictly decreasing in $x \geq 0$, then $\criticalpoint{\eta}$
is unique.
\end{proposition}

Recall that the optimality of a threshold policy 
should not come as a surprise, and could in fact be shown in the pre-limit
and within a far wider class of policies than those satisfying \refEquation{formf}.
The strength of \refProposition{prop:Threshold_control_is_optimal_form_of_control} lies in the characterization \refEquation{eqn:Threshold_equation} of
$\criticalpoint{\eta}$. We refer to \refEquation{eqn:Threshold_equation} as the \emph{threshold equation}: it is a powerful basis on which to obtain numerical
solutions, closed-form expressions, bounds, and asymptotic expansions for
$\criticalpoint{\eta}$. Results for $\criticalpoint{\eta}$ of this nature are presented in \refSection{sec:Analysis_of_optimal_thresholds}.

\begin{example}[Exponential revenue revisited]
\label{eg:Exponential_revenue_revisited}
Let us revisit \refExample{expo}, where $r(x) = \exp{(bx)}$ for $x < 0$
and $r(x) = \exp{(-dx)}$ for $x \geq 0$. The threshold equation, \refEquation{eqn:Threshold_equation}, takes the form
\begin{equation}
\e{-d \eta} \Bigl( \frac{\Phi(\gamma)}{\phi(\gamma)} + \frac{1 - \e{-\gamma \eta}}{\gamma} \Bigr)
= \frac{\Phi(\gamma-b)}{\phi(\gamma-b)} + \frac{ 1 - \e{-(d+\gamma)\eta} }{d+\gamma},
\label{eqn:Example__Exponential_revenue__Threshold_equation}
\end{equation}
which we study in depth in \refSection{sec:Exact_characterizations__Exponential_revenue}.
When $b = 5$, $d=1$, and $\gamma = 0.01$, solving \refEquation{eqn:Example__Exponential_revenue__Threshold_equation} numerically yields $\criticalpoint{\eta} \approx 1.00985$, which supports our earlier observation that $\criticalpoint{\eta} \approx 1$ in \refExample{expo}.
\end{example}

In \refFigure{fig:Figure__Optimal_buffer_size_for_finite_systems_as_s_increases}, the true optimal admission threshold
$\criticalpoint{\tau} = \arg\max_{\tau \in \naturalNumbersZero} R_{\mathrm{T},s}(\tau)$ is plotted as a function of $s$, along with the asymptotic \gls{QED} approximation 
$\criticalpoint{\tau} \approx \lfloor \criticalpoint{\eta} \sqrt{s} \rfloor$.
We observe that the \gls{QED} approximation is accurate,
even for a relatively small number of servers, and exact in the
majority of cases. This is reflected in \refFigure{fig:Figure__Optimal_revenue_for_finite_systems_as_s_increases__Relative_optimality_gap},
which plots the relative optimality gap as a function of $s$.
The relative optimality gap is zero for the vast majority of $s$ values,
and as low as $10^{-2}$ for systems with as few as $10$ servers.

\begin{figure}[!hbtp]
\begin{center}
\small
\input{Figure__Optimal_buffer_size_for_finite_systems_as_s_increases}
\end{center}
\vspace{-1em}
\caption{\textrm{The true optimal admission threshold $\criticalpoint{\tau}$
as a function of $s$, together with the (almost indistinguishable) \gls{QED}
approximation $\lfloor \criticalpoint{\eta} \sqrt{s} \rfloor$.}}
\label{fig:Figure__Optimal_buffer_size_for_finite_systems_as_s_increases}
\end{figure}

\begin{figure}[!hbtp]
\begin{center}
\small
%\documentclass{standalone}
%
%\usepackage{amsmath}
%\usepackage{amsthm}
%\usepackage{amssymb}
%\usepackage{siunitx}
%
%\usepackage{tikz}
%\usepackage{pgfplots}
%\usepgfplotslibrary{groupplots}
%
%\newcommand{\criticalpoint}[1]{  #1^{\textnormal{opt}} }
%
%\begin{document}
\begin{tikzpicture}
\begin{axis}[
	width=0.9\columnwidth, height = 0.5*0.9*0.618\columnwidth,
%	width=\columnwidth, height=0.618\columnwidth,
    xmin=0.85, xmax=275,
    xmode=log,
    ymode=log,
    ymin=0.0000005, ymax=2,
    xtick={1,10,100},
    every axis x label/.style=
        {at={(ticklabel cs:0.5)},anchor=north},
    ytick={0,0.000001,0.0001,0.01,1},
    y tick label style={/pgf/number format/fixed},
    every axis y label/.style=
%        {at={(ticklabel cs:0.5)},rotate=90,anchor=south},
        {at={(ticklabel cs:0.5)},rotate=0,anchor=east},
    scaled ticks=true,
	xlabel={$s$},  			
	]

%%%%%%%%%%%%%%%%%%%%%%%%%%%
%% The threshold graphs. %%
%%%%%%%%%%%%%%%%%%%%%%%%%%%

\addplot[color=black, draw=none, thick, mark=+] plot coordinates {
(1, 0.08762) (2, 0.) (3, 0.) (4, 0.0149866) (5, 0.) (6, 0.) (7, 0.) (8, 0.) (9, 0.00365697) (10, 0.) (11, 0.) (12, 0.) (13, 0.) (14, 0.) (15, 0.) (16, 0.000865061) (17, 0.) (18, 0.) (19, 0.) (20, 0.) (21, 0.) (22, 0.) (23, 0.) (24, 0.) (25, 0.) (26, 0.) (27, 0.) (28, 0.) (29, 0.) (30, 0.) (31, 0.) (32, 0.) (33, 0.) (34, 0.) (35, 0.) (36, 0.) (37, 0.) (38, 0.) (39, 0.) (40, 0.) (41, 0.) (42, 0.) (43, 0.) (44, 0.) (45, 0.) (46, 0.) (47, 0.) (48, 0.) (49, 0.) (50, 0.) (51, 0.) (52, 0.) (53, 0.) (54, 0.) (55, 0.) (56, 0.) (57, 0.) (58, 0.) (59, 0.) (60, 0.) (61, 0.) (62, 0.) (63, 0.) (64, 0.) (65, 0.) (66, 0.) (67, 0.) (68, 0.) (69, 0.) (70, 0.) (71, 0.) (72, 0.) (73, 0.) (74, 0.) (75, 0.) (76, 0.) (77, 0.) (78, 0.) (79, 0.) (80, 0.) (81, 0.) (82, 0.) (83, 0.) (84, 0.) (85, 0.) (86, 0.) (87, 0.) (88, 0.) (89, 0.) (90, 0.) (91, 0.) (92, 0.) (93, 0.) (94, 0.) (95, 0.) (96, 0.) (97, 0.) (98, 0.) (99, 0.) (100, 0.) (101, 0.) (102, 0.) (103, 0.) (104, 0.) (105, 0.) (106, 0.) (107, 0.) (108, 0.) (109, 0.) (110, 0.) (111, 0.) (112, 0.) (113, 0.) (114, 0.) (115, 0.) (116, 0.) (117, 0.) (118, 0.) (119, 0.) (120, 0.) (121, 0.) (122, 0.) (123, 0.) (124, 0.) (125, 0.) (126, 0.) (127, 0.) (128, 0.) (129, 0.) (130, 0.) (131, 0.) (132, 0.) (133, 0.) (134, 0.) (135, 0.) (136, 0.) (137, 0.) (138, 0.) (139, 0.) (140, 0.) (141, 0.000000876075) (142, 0.) (143, 0.) (144, 0.) (145, 0.) (146, 0.) (147, 0.) (148, 0.) (149, 0.) (150, 0.) (151, 0.) (152, 0.) (153, 0.) (154, 0.) (155, 0.) (156, 0.) (157, 0.) (158, 0.) (159, 0.) (160, 0.) (161, 0.) (162, 0.) (163, 0.) (164, 0.) (165, 0.) (166, 0.) (167, 0.) (168, 0.) (169, 0.) (170, 0.) (171, 0.) (172, 0.) (173, 0.) (174, 0.) (175, 0.) (176, 0.) (177, 0.) (178, 0.) (179, 0.) (180, 0.) (181, 0.) (182, 0.) (183, 0.) (184, 0.) (185, 0.) (186, 0.) (187, 0.) (188, 0.) (189, 0.) (190, 0.) (191, 0.) (192, 0.0000183047) (193, 0.) (194, 0.) (195, 0.) (196, 0.) (197, 0.) (198, 0.) (199, 0.) (200, 0.) (201, 0.) (202, 0.) (203, 0.) (204, 0.) (205, 0.) (206, 0.) (207, 0.) (208, 0.) (209, 0.) (210, 0.) (211, 0.) (212, 0.) (213, 0.) (214, 0.) (215, 0.) (216, 0.) (217, 0.) (218, 0.) (219, 0.) (220, 0.) (221, 0.) (222, 0.) (223, 0.) (224, 0.) (225, 0.) (226, 0.) (227, 0.) (228, 0.) (229, 0.) (230, 0.) (231, 0.) (232, 0.) (233, 0.) (234, 0.) (235, 0.) (236, 0.) (237, 0.) (238, 0.) (239, 0.) (240, 0.) (241, 0.) (242, 0.) (243, 0.) (244, 0.) (245, 0.) (246, 0.) (247, 0.) (248, 0.) (249, 0.) (250, 0.) (251, 0.0000327945) (252, 0.) (253, 0.) (254, 0.) (255, 0.) (256, 0.) 
};
\draw (4,-5) node [fill=none, draw=none, text=black, opacity=1, text opacity=1, scale=1.25] {$\frac{ R_{\mathrm{T},s}(\criticalpoint{\tau}) - R_{\mathrm{T},s}( \lfloor \criticalpoint{\eta} \sqrt{s} \rfloor ) }{ R_{\mathrm{T},s}(\criticalpoint{\tau}) }$};

\end{axis}
\end{tikzpicture}
%\end{document}
\end{center}
\vspace{-1em}
\caption{\textrm{The relative error $( R_{\mathrm{T},s}(\criticalpoint{\tau}) - R_{\mathrm{T},s}( \lfloor \criticalpoint{\eta} \sqrt{s} \rfloor ) )/ R_{\mathrm{T},s}(\criticalpoint{\tau})$ as a function of $s$. The missing points indicate an error that is strictly zero. The errors that are non-zero arise due to the \gls{QED} approximation for the optimal
admission threshold being off by just one state.}}
\label{fig:Figure__Optimal_revenue_for_finite_systems_as_s_increases__Relative_optimality_gap}
\end{figure}

We remark that when utilizing the asymptotic optimal threshold provided by \refProposition{prop:Threshold_control_is_optimal_form_of_control} in a finite system, the proof of \refProposition{prop:QED_limit_of_revenue__Continuous_admission_controls} in \refAppendixSection{sec:Limiting_behavior_of_long_term_QED_revenue} guarantees that $R_{s,\mathrm{T}}( \lfloor \criticalpoint{\eta} \sqrt{s} \rfloor ) - R_{\mathrm{T}}( \criticalpoint{\eta} ) = \bigO{ 1 / \sqrt{s} }$. In other words, a finite system that utilizes the asymptotic optimal threshold achieves a revenue within $\bigO{1/\sqrt{s}}$ of the solution to \refEquation{eqn:Objective_optimization_problem_postlimit}.

\subsection{Customer reward maximization}
\label{sec:Revenue_framework__Customer_state_dependent_rewards}

In \refSection{sec:Introduction} we have discussed the difference between revenues seen from the system's perspective and from the customer's perspective. Although the emphasis lies on the system's perspective, as in \refSection{sec:Revenue_framework__Asymptotic_control_problem}, we now show how results for the customer's perspective can be obtained.

\subsubsection{Linear revenue structure}

If revenue is generated at rate $a > 0$ for each customer that is being served, and cost is incurred at rate $b > 0$ for
each customer that is waiting for service, the revenue structure is given by
\begin{equation}
r_s(k)
=
\begin{cases}
a k, & k \leq s, \\
a s - b (k - s), & k \geq s.
\end{cases}
\label{linestru}
\end{equation}
When $n_s = a s$ and $q_s = \sqrt{s}$, the revenue structure in \refEquation{linestru} satisfies the scaling condition in \refEquation{formr}, with
\begin{equation}
r(x) =
\begin{cases}
a x, & x \leq 0, \\
- b x, & x \geq 0. \\
\end{cases}
\label{lineprof}
\end{equation}
Consequently, \refProposition{prop:QED_limit_of_revenue__Continuous_admission_controls} implies that
\begin{equation}
\lim_{s \rightarrow \infty} \frac{ R_s(\{p_s(k)\}_{k\geq0}) - a s }{\sqrt{s}} 
%= \frac{ - a \int_0^\infty x \e{- \frac{1}{2} x^2 + \gamma x} \d{x} - b \int_0^\infty x f(x) \e{- \gamma x} \d{x}}{ \frac{\Phi(\gamma)}{\phi(\gamma)} + \int_0^\infty f(x) \e{- \gamma x} \d{x} }
= \frac{ a \bigl( 1 + \gamma \frac{\Phi(\gamma)}{\phi(\gamma)} \bigr) -
b \int_0^\infty x f(x) \e{- \gamma x} \d{x}}
{\frac{\Phi(\gamma)}{\phi(\gamma)} + \int_0^\infty f(x) \e{- \gamma x} \d{x}},
\end{equation}
for any profile $f$, and \refProposition{prop:Threshold_control_is_optimal_form_of_control} reveals that
the optimal control is $f(x) = \indicator{ 0 \leq x < \criticalpoint{\eta} }$ with $\criticalpoint{\eta}$ the unique solution of the threshold equation \refEquation{eqn:Threshold_equation}, which with $c=a/b$ becomes
\begin{equation}
\eta \Bigl( \frac{\Phi(\gamma)}{\phi(\gamma)} + \frac{ 1 - \e{-\gamma \eta} }{\gamma} \Bigr) 
= c \Bigl( 1 + \gamma \frac{\Phi(\gamma)}{\phi(\gamma)} \Bigr) + \frac{ 1- (1+\gamma \eta) \e{-\gamma \eta} }{ \gamma^2 }.
\label{eqn:Example__Linear_revenue__Threshold_equation}
\end{equation}
We see that $\criticalpoint{\eta}$ depends only on $c$. The threshold equation \refEquation{eqn:Example__Linear_revenue__Threshold_equation} is studied extensively in \refSection{sec:Exact_characterizations__Linear_revenue}. A minor variation of the arguments used there to prove
\refProposition{prop:Optimal_threshold_for_linear_revenue}, shows that
\begin{equation}
\criticalpoint{\eta} 
= r_0 + \frac{1}{\gamma} W\Bigl( \frac{\gamma \e{-\gamma r_0}}{a_0} \Bigr),
\label{etaopt}
\end{equation} 
where $W$ denotes the Lambert W function, and
\begin{equation}
a_0 
%= - \gamma^2 \Bigl(\frac{\Phi(\gamma)}{\phi(\gamma)} + \frac{1}{\gamma^2}\Bigr) 
= - 1 - \gamma^2  \frac{\Phi(\gamma)}{\phi(\gamma)},
\quad
r_0 
%= \frac{c \bigl(1 + \gamma \frac{\Phi(\gamma)}{\phi(\gamma)}\bigr) + \frac{1}{\gamma^2}}{\frac{\Phi(\gamma)}{\phi(\gamma)} + \frac{1}{\gamma}} 
= c \gamma + \frac{1}{\gamma + \gamma^2 \frac{\Phi(\gamma)}{\phi(\gamma)}}.
\end{equation}

\subsubsection{Relating system-governed revenue to customer rewards}

From \refEquation{etaopt}, it can be deduced that for large values of $\gamma$,
$\criticalpoint{\eta} \approx c \gamma$ (see the proof of \refProposition{prop:Asymptotics_of_LambertW_solutions_as_gamma_to_infinity}, and for a discussion on the asymptotic behavior of the threshold equation for general revenue structures, we refer to \refSection{sec:Analyses__Asymptotic_solutions}). Thus, asymptotically, the optimal threshold value is approximately equal to the product of
the staffing slack $\gamma$ and the ratio of the service revenue $a$
and the waiting cost $b$.
%For a discussion on the asymptotic behavior of the threshold equation for general revenue structures, we refer the reader to \refSection{sec:Analyses__Asymptotic_solutions}. 
%The asymptotic behavior of the optimal threshold value $\criticalpoint{\eta} \approx c \gamma$
%for large values of $\gamma$ in the case of linear revenue structures may be explained as follows.

The asymptotic behavior $\criticalpoint{\eta} \approx c \gamma$ may be explained as follows.
For each arriving customer, we must balance the expected revenue $a$
when that customer is admitted and eventually served against the expected
waiting cost incurred for that customer as well as the additional waiting cost
for customers arriving after that customer.
When the arriving customer finds $\tau$ customers waiting, the overall waiting
cost incurred by admitting that customer may be shown to behave roughly as
$b \tau / (\gamma \sqrt{s})$ for large values of $\gamma$.
Equating $a$ with $b \tau / (\gamma \sqrt{s})$ then yields that the optimal
threshold value should approximately be $\criticalpoint{\tau} \approx c \gamma \sqrt{s}$.

As the above explanation reveals, the stationary average revenue rate
$R_s(\{p_s(k)\}_{k\geq0})$ under revenue structure \refEquation{linestru} is the same
as when a reward $a > 0$ is received for each admitted customer and a penalty
$b \expectation{W}$ is charged when the expected waiting time of that customer is $\expectation{W}$, with $b > 0$.
In the latter case the stationary average reward earned may be expressed as \refEquation{eqn:Customer_oriented_reward_rate} in, where now
\begin{equation}
\hat{r}_s(k) = a - b \max \Bigl\{0, \frac{k - s + 1}{s} \Bigr\}
\label{reward}
\end{equation}
denotes a \emph{customer reward}.

The \emph{system-governed} revenue rate and the
\emph{customer} reward rate are in this case equivalent. To see this, write
\begin{equation}
\hat{R}_s(\{p_s(k)\}_{k\geq0})
= a \lambda \Bigl( \sum_{k = 0}^{s - 1} \pi_s(k) + \sum_{k = s}^{\infty} p_s(k-s) \pi_s(k) \Bigr) - b \lambda \sum_{k = s}^{\infty} \frac{k - s + 1}{s} p_s(k-s) \pi_s(k).
\end{equation}
Then note that because the arrival rate multiplied by the probability that an arriving customer is admitted must equal the expected number of busy servers, and by local balance $\lambda \pi_s(k) p_s(k-s) = s \pi_s(k + 1)$ for $k = s, s + 1, \dots$, we have
\begin{align}
&
\hat{R}_s(\{p_s(k)\}_{k\geq0})
= a \Bigl( \sum_{k = 0}^{s - 1} k \pi_s(k) + \sum_{k = s}^{\infty} s \pi_s(k) \Bigr) - b \sum_{k = s + 1}^{\infty} (k - s) \pi_s(k) 
\nonumber \\ &
= \sum_{k = 0}^{s - 1} a k \pi_s(k) + \sum_{k = s}^{\infty} ( a s - b (k - s) ) \pi_s(k)
\overset{\refEquation{linestru}} = \sum_{k = 0}^{\infty} r_s(k) \pi_s(k) 
= R_s(\{p_s(k)\}_{k\geq0}).
\end{align}
The optimal threshold in \refEquation{etaopt} thus maximizes the customer reward rate $\hat{R}_s(\{p_s(k)\}_{k\geq0})$ asymptotically as well, i.e., in this example by \refProposition{prop:Threshold_control_is_optimal_form_of_control}, 
\begin{equation}
\lim_{s \to \infty} \frac{ \hat{R}_s(\{p_s(k)\}_{k\geq0}) - a }{ \sqrt{s} }
= \lim_{s \to \infty} \frac{ R_s(\{p_s(k)\}_{k\geq0}) - a }{ \sqrt{s} }
= R(f).
\end{equation}

In fact, for \emph{any} system-governed revenue rate $r_s(k)$, the related customer reward structure $\hat{r}_s(k)$ is given by
\begin{equation}
\hat{r}_s(k) 
= \frac{ r_s(k+1) }{ \min\{k+1,s\} },
%= \frac{r_s(k+1)}{ \# \textrm{ busy servers in state } k+1 },
\quad k \in \naturalNumbersZero,
\label{eqn:Relation_between_customer_oriented_reward_structures_and_system_governed_revenue}
\end{equation}
because then
\begin{align}
&
\hat{R}_s(\{p_s(k)\}_{k\geq0})
= \sum_{k=0}^{s-1} \hat{r}_s(k) \lambda \pi_s(k) + \sum_{k=s}^\infty \hat{r}_s(k) \lambda p_s(k-s) \pi_s(k)
%\nonumber \\ &
= \sum_{k=0}^{s-1} \hat{r}_s(k) (k+1) \pi_s(k+1) 
\nonumber \\ &
+ \sum_{k=s}^\infty \hat{r}_s(k) s \pi_s(k+1)
= \sum_{k=0}^{s} \hat{r}_s(k-1) k \pi_s(k) + \sum_{k=s+1}^\infty \hat{r}_s(k-1) s \pi_s(k)
%\nonumber \\ &
%= \sum_{k=0}^\infty r_s(k) \pi_s(k) 
%\nonumber \\ &
= R_s(\{p_s(k)\}_{k\geq0}),
\end{align}
using local balance, $\lambda \pi_s(k) = (k+1) \pi_s(k+1)$ for $k = 0, 1, \ldots, s-1$ and $\lambda p_s(k-s) \pi_s(k) = s \pi_s(k+1)$ for $k = s, s+1, \ldots$.

\begin{proposition}
\label{prop:Relation_between_system_governed_revenue_and_customer_rewards}
For any system-governed revenue rate $r_s(k)$, the customer reward structure $\hat{r}_s(k)$ in \refEquation{eqn:Relation_between_customer_oriented_reward_structures_and_system_governed_revenue} guarantees that the average system-governed revenue rate $R_s(\{p_s(k)\}_{k\geq0})$ equals the customer reward rate $\hat{R}_s(\{p_s(k)\}_{k\geq0})$.
\end{proposition}

In particular, \refProposition{prop:Relation_between_system_governed_revenue_and_customer_rewards} implies that counterparts to \refProposition{prop:QED_limit_of_revenue__Continuous_admission_controls} and \refProposition{prop:Threshold_control_is_optimal_form_of_control} hold for the customer reward rate $\hat{R}_s(\{p_s(k)\}_{k\geq0})$, assuming the customer reward structure $\hat{r}_s(k)$ is appropriately scaled.

\section{Properties of the optimal threshold}
\label{sec:Analysis_of_optimal_thresholds}

%Because customer reward structures can be related to system-governed revenue rates through \refEquation{eqn:Relation_between_customer_oriented_reward_structures_and_system_governed_revenue}, we may focus on maximization of either of them.

We focus throughout this paper on maximization of the average system-governed revenue rate. In \refSection{sec:Revenue_framework} we have established threshold optimality and derived the threshold equation that defines the optimal threshold $\criticalpoint{\eta}$. In this section we obtain a series of results about $\criticalpoint{\eta}$. In \refSection{sec:Analyses__General_bounds}, we present a procedure (for general revenue functions) to obtain an upper bound $\eta^{\max}$ and a lower bound $\eta^{\min}$ on the optimal threshold $\criticalpoint{\eta}$. \refSection{sec:Analyses__Monotonicity} discusses our monotonicity results. Specifically, we prove that $\criticalpoint{\eta}$ increases with $\gamma \in \realNumbers$, and that $R_{\mathrm{T}}(0)$, $R_{\mathrm{T}}(\criticalpoint{\eta})$ both decrease with $\gamma \in \realNumbers$. In \refSection{sec:Analyses__Asymptotic_solutions}, we derive asymptotic descriptions of the optimal threshold for general revenue structures, even if the revenue structures would not allow for an explicit characterization. We prove that $\criticalpoint{\eta} \approx r^{\gets}( r(-\gamma) )$ as $\gamma \rightarrow \infty$, and that $\criticalpoint{\eta} \approx -(1/\gamma) \ln{( 1 - r'(0-) / r'(0+) )}$ as $\gamma \rightarrow - \infty$. In \refSection{sec:Analyses__Exact_characterizations}, we derive explicit characterizations of $\criticalpoint{\eta}$ for linear and exponential revenue structures.

%\subsection{Preliminaries}

We henceforth assume that $r(x)$ is piecewise smooth and bounded on $(-\infty,0)$ and $(0,\infty)$, and continuous at $0$ with
$r(\pm 0) = 1 = r(0)$. We also assume that $r(x)$ is increasing on $(-\infty,0]$ and decreasing on $[0,\infty)$, with $0 \leq r(x) \leq r(0) = 1$. Revenue functions for which $r(0) > 0$ and $r(0) \neq 1$ can be considered through analysis of the scaled revenue function $\bar{r}(x) = r(x) / r(0)$. For notational convenience, we also define $r_{\mathrm{L}}(x)$ and $r_{\mathrm{R}}(x)$ as
\begin{equation}
r(x)
=
\begin{cases}
r_{\mathrm{L}}(x), & x < 0, \\
r_{\mathrm{R}}(x), & x \geq 0, \\
\end{cases}
\end{equation}
and introduce $A = \int_{-\infty}^0 r_{\mathrm{L}}(x) \e{-\frac{1}{2}x^2-\gamma x} \d{x}$, and $B = \Phi(\gamma) / \phi(\gamma)$.  Note that $R_{\mathrm{T}}(0) = A / B$.

\begin{corollary}
Under these assumptions, there exists a solution $\criticalpoint{\eta} > 0$ of the threshold equation. This solution is positive, unless $r_{\mathrm{L}}(x) = 1$, and unique when $r_{\mathrm{R}}'(x) < 0$ for all $x \geq 0$ such that $r_{\mathrm{R}}(x) > 0$.
\end{corollary}

\myProof{
Note that these assumptions on $r$ are slightly stronger than in \refProposition{prop:Threshold_control_is_optimal_form_of_control}. This corollary is directly implied by our proof of \refProposition{prop:Threshold_control_is_optimal_form_of_control}; see the explanation between \refEquation{eqn:Definition__Set_x0}
and \refEquation{eqn:Bounds__On_integral_over_sx_gx}.
}

\subsection{General bounds}
\label{sec:Analyses__General_bounds}

Denote the inverse function of $r$ by $r^{\gets}$. The following bounds hold for general revenue structures, are readily calculated, and provide insight into the optimal thresholds. Later in \refSection{sec:Exact_characterizations__Linear_revenue} and \refSection{sec:Exact_characterizations__Exponential_revenue}, we illustrate \refProposition{prop:Threshold_control_bounds} for a linear revenue structure, and an exponential revenue structure, respectively.

\begin{proposition}
\label{prop:Threshold_control_bounds}
When $r$ is strictly decreasing for $x \geq 0$, $\eta^{\max} = r_{\mathrm{R}}^{\gets}( R_{\mathrm{lower}} ) \geq \criticalpoint{\eta}$, and $\eta^{\min} = r_{\mathrm{R}}^{\gets}( R_{\mathrm{upper}} ) < \criticalpoint{\eta}$. Here,
\begin{equation}
R_{\mathrm{lower}} 
= R_{\mathrm{T}}(0) 
= \frac{ \int_{-\infty}^0 r(x) \e{ - \frac{1}{2} x^2 - \gamma x } \d{x} }{ \frac{\Phi(\gamma)}{\phi(\gamma)} },
\quad 
R_{\mathrm{upper}} 
= \frac{ \int_{-\infty}^0 r(x) \e{ - \frac{1}{2} x^2 - \gamma x } \d{x} + \int_0^{\eta^{\max}} \e{-\gamma x} \d{x} }{ \frac{\Phi(\gamma)}{\phi(\gamma)} + \int_0^{\eta^{\max}} \e{-\gamma x} \d{x} }.
\end{equation}
\end{proposition}

\myProof{
The assumptions on $r_{\mathrm{R}}(x)$ imply that its inverse function $r_{\mathrm{R}}^{\gets}(y)$ exists, and that it is also strictly decreasing. It is therefore sufficient to provide upper and lower bounds on $R_{\mathrm{T}}(\eta)$ that are independent of $\eta$.

For threshold control policies, the system revenue is given by \refEquation{eqn:Long_term_revenue__Threshold}. Also recall that the optimal threshold $\criticalpoint{\eta}$ solves the threshold equation, i.e., $r_{\mathrm{R}}( \criticalpoint{\eta} ) = R_{\mathrm{T}}( \criticalpoint{\eta} )$. By suboptimality, we immediately obtain $R_{\mathrm{lower}} \leq R_{\mathrm{T}}(\criticalpoint{\eta})$, and so $\criticalpoint{\eta} \leq \eta^{\max}$ by monotonicity. 

We will first derive an alternative forms of the threshold equation. For instance, rewriting \refEquation{eqn:Threshold_equation} into
\begin{equation}
\Bigl( B + \int_0^\eta \e{-\gamma x} \d{x} \Bigr) r(\eta) 
= A + \int_{0}^\eta r(x) \e{ - \gamma x } \d{x},
\label{eqn:Threshold_equation_rewritten}
\end{equation}
dividing by $B$, and using $R_{\mathrm{T}}(0) = A / B$ gives
\begin{equation}
r(\eta) - R_{\mathrm{T}}(0) 
= - \frac{r(\eta)}{ B } \int_0^\eta \e{-\gamma x} \d{x} + \int_{0}^\eta r(x) \e{ - \gamma x } \d{x}.
\end{equation}
We then identify the right-hand member as being a result of an integration by parts, to arrive at the alternative form
\begin{equation}
r(\eta) - R_{\mathrm{T}}(0) 
= - \Bigl[ \frac{r(x)}{ B } \int_0^x \e{-\gamma u} \d{u} \Big]_0^\eta + \int_0^\eta \e{-\gamma x} \d{x}
= - \int_{0}^\eta \frac{r'(x)}{ B } \int_0^x \e{-\gamma u} \d{u} \d{x}.
\label{eqn:Alternative_form_of_the_threshold_equation}
\end{equation}

Let $c(\eta) = \int_0^\eta \e{-\gamma x} \d{x} = (1-\e{-\gamma\eta})/\gamma$ if $\gamma \neq 0$ and $c(\eta) = \eta$ if $\gamma = 0$. Since $c(\eta)$ is increasing in $\eta$ and $-r'(x) \geq 0$ for $x \geq 0$, we have for $\eta \geq 0$
\begin{equation}
- \frac{1}{B} \int_0^\eta r'(x) c(x) \d{x}
< - \frac{1}{B} c(\eta) \int_0^\eta r'(x) \d{x}
= \frac{1}{B} c(\eta) ( 1 - r(\eta) ).
\end{equation}
Let $\eta = \hat{\eta}$ be the (unique) solution of the equation 
\begin{equation}
r(\eta) - R_{\mathrm{T}}(0) 
= \frac{1}{B} c(\eta) ( 1 - r(\eta) ).
\label{eqn:Equation_for_eta_hat}
\end{equation}
Then 
\begin{equation}
r(\hat{\eta}) - R_{\mathrm{T}}(0)
= \frac{1}{B} c(\hat{\eta}) ( 1 - r(\hat{\eta}) ) 
> - \frac{1}{B} \int_0^{\hat{\eta}} r'(x) c(x) \d{x},
\end{equation}
and so $0 < \hat{\eta} < \criticalpoint{\eta}$. We have from \refEquation{eqn:Equation_for_eta_hat} that
\begin{equation}
r(\hat{\eta})
= \frac{ R_{\mathrm{T}}(0) + \frac{1}{B} c(\hat{\eta}) }{ 1 + \frac{1}{B} c(\hat{\eta}) }
= \frac{ A + c(\hat{\eta}) }{ B + c(\hat{\eta}) }.
\label{eqn:Expression_for_r_evaluated_at_eta_hat}
\end{equation}
From $\hat{\eta} < \criticalpoint{\eta} < \eta^{\max} = r^{\gets}(R_{\mathrm{T}}(0))$, we then find
\begin{equation}
c(\hat{\eta}) < c(\eta^{\max}),
\quad
\frac{ A + c(\hat{\eta}) }{ B + c(\hat{\eta}) }
< \frac{ A + c(\eta^{\max}) }{ B + c(\eta^{\max}) }
= R_{\mathrm{upper}},
\end{equation}
since $0 < A < B$, i.e.\ $R_{\mathrm{T}}(0) < 1$. This completes the proof.
}

\subsection{Monotonicity}
\label{sec:Analyses__Monotonicity}

We next investigate the influence of the slack $\gamma$ on the optimal threshold.

\begin{proposition}
\label{prop:Monotonicity_of_R0}
The revenue $R_{\mathrm{T}}(0)$ decreases in $\gamma \in \realNumbers$.
\end{proposition}

\myProof{
Write $r(-x) = u(x)$ so that $0 \leq u(x) \leq 1 = u(0)$ and
\begin{equation}
R_{\mathrm{T}}(0) 
= \frac{ \int_0^\infty u(x) \e{-\frac{1}{2}x^2+\gamma x} \d{x} }{ \int_0^\infty \e{-\frac{1}{2}x^2+\gamma x} \d{x} },
\quad \gamma \in \realNumbers,
\end{equation}
and calculate
\begin{equation}
\frac{ \d{R_{\mathrm{T}}(0)} }{ \d{\gamma} }  
= \frac{ \int_0^\infty \e{-\frac{1}{2}x^2+\gamma x} \d{x} \int_0^\infty x u(x) \e{-\frac{1}{2}x^2+\gamma x} \d{x} - \int_0^\infty u(x) \e{-\frac{1}{2}x^2+\gamma x} \d{x} \int_0^\infty x \e{-\frac{1}{2}x^2+\gamma x} \d{x} }{ \Bigl( \int_0^\infty \e{-\frac{1}{2}x^2+\gamma x} \d{x} \Bigr)^2 }.
\end{equation}
The numerator can be written as
\begin{equation}
N 
= \int_0^\infty \int_0^\infty (x-y) u(x) \e{-\frac{1}{2}x^2+\gamma x} \e{-\frac{1}{2}y^2+\gamma y} \d{x} \d{y}.
\end{equation}
Suppose that $u(x)  = \indicator{ 0 \leq x < a }$ for some $a > 0$. Then
\begin{equation}
N 
= \int_0^a \int_0^\infty (x-y) \e{-\frac{1}{2}x^2+\gamma x} \e{-\frac{1}{2}y^2+\gamma y} \d{x} \d{y} 
\leq \int_0^a \int_0^a (x-y) \e{-\frac{1}{2}x^2+\gamma x} \e{-\frac{1}{2}y^2+\gamma y} \d{x} \d{y} 
= 0.
\end{equation}
In general, we can write $u(x) = -\int_0^\infty \indicator{ 0 \leq x < a } u'(a) \d{a} = - \int_x^\infty u'(a) \d{a}$ with $u'(a) < 0$, to arrive at
\begin{equation}
N 
= \int_0^\infty - u'(a) \Bigl( \int_0^\infty \int_0^\infty (x-y) \indicator{ 0 \leq x < a } \e{-\frac{1}{2}x^2+\gamma x} \e{-\frac{1}{2}y^2+\gamma y} \d{x} \d{y} \Bigr) \d{a} 
\leq 0.
\end{equation}
This concludes the proof.
}

\begin{proposition}
\label{prop:Optimal_threshold_is_increasing_in_gamma}
The optimal threshold $\criticalpoint{\eta}$ increases in $\gamma \in \realNumbers$, and $R_{\mathrm{T}}(\criticalpoint{\eta})$ decreases in $\gamma \in \realNumbers$.
\end{proposition}

\myProof{
By \refProposition{prop:Monotonicity_of_R0}, we have that $R_{\mathrm{T}}(0)$ decreases in $\gamma \in \realNumbers$. Furthermore, for any $\eta > 0$, we have that
$
\int_0^\eta \e{-\gamma x} \d{x} / B
$
decreases in $\gamma \in \realNumbers$. Consider the alternative form \refEquation{eqn:Alternative_form_of_the_threshold_equation} of the threshold equation. For fixed $\eta$, the left member thus increases in $\gamma$, while the right member decreases in $\gamma$, since $r'(x) < 0$. The solution $\criticalpoint{\eta}$ of the threshold equation therefore increases in $\gamma \in \realNumbers$.

To prove the second part of the claim, we recall that $\criticalpoint{\eta}$ solves the threshold equation, so $R_{\mathrm{T}}(\criticalpoint{\eta}) = r_{\mathrm{R}}(\criticalpoint{\eta})$. Since $\criticalpoint{\eta} \geq 0$ is increasing in $\gamma \in \realNumbers$, our assumptions on $r$ imply that $r_{\mathrm{R}}(\criticalpoint{\eta})$ is decreasing in $\gamma \in \realNumbers$. Hence, $R_{\mathrm{T}}(\criticalpoint{\eta})$ is decreasing in $\gamma \in \realNumbers$ as well.
}

\refProposition{prop:Optimal_threshold_is_increasing_in_gamma} can be interpreted as follows. First note that an increase in $\gamma$ means that fewer customers are served by the system, apart from the impact of a possible admission control policy. Then, for threshold control, an increased $\gamma$ implies that the optimal threshold should increase, in order to serve more customers. This of course is a direct consequence of our revenue structure, which is designed to let the system operate close to the ideal operating point. A large $\gamma$ drifts the process away from this ideal operating point, and this can be compensated for by a large threshold $\criticalpoint{\eta}$. Hence, although the slack $\gamma$ and the threshold $\criticalpoint{\eta}$ have quite different impacts on the system behavior, at a high level their monotonic relation can be understood, and underlines that the revenue structure introduced in this paper has the right properties for the \gls{QED} regime.

\subsection{Asymptotic solutions}
\label{sec:Analyses__Asymptotic_solutions}

We now present asymptotic results for the optimal threshold in the regimes where the slack $\gamma$ becomes extremely large or extremely small.

\begin{proposition}
\label{prop:Asymptotic_solutions}
When $\gamma \rightarrow -\infty$, and if the revenue function has a cusp at $x = 0$, i.e., $r_{\mathrm{R}}'(0+) < 0 < r_{\mathrm{L}}'(0-)$, the optimal threshold is given by
\begin{equation}
\criticalpoint{\eta}
= - \frac{1}{\gamma} \ln{ \Bigl( 1 - \frac{r_{\mathrm{L}}'(0-)}{r_{\mathrm{R}}'(0+)} \Bigr) } + \bigObig{ \frac{1}{\gamma^2} }.
\label{eqn:Extreme_gamma__Threshold_equation_when_gamma_is_near_minus_infinity}
\end{equation}
\end{proposition}

\myProof{
We consider $\gamma \to - \infty$. From steepest descent analysis, we have for a smooth and bounded $f$ on $(-\infty,0]$,
\begin{equation}
\int_{-\infty}^0 f(x) \e{-\gamma x} \d{x} 
= - \frac{f(0)}{\gamma} - \frac{f'(0)}{\gamma^2} + \bigObig{ \frac{1}{\gamma^3} }, 
\quad \gamma \to -\infty.
\end{equation}
Hence, it follows that
\begin{equation}
R_{\mathrm{T}}(0)
= \frac{A}{B}
= \frac{ -\frac{1}{\gamma} - \frac{r_{\mathrm{L}}'(0-)}{\gamma^2} + \bigObig{ \frac{1}{\gamma^3} } }{ -\frac{1}{\gamma} + \bigObig{ \frac{1}{\gamma^3} } }
= 1 + \frac{r_{\mathrm{L}}'(0-)}{\gamma} + \bigObig{ \frac{1}{\gamma^2} },
\quad \gamma \to -\infty.
\label{eqn:Expansion_for_RT0_as_gamma_to_minus_infty}
\end{equation}
From the upper bound $\criticalpoint{\eta} < r^{\gets}(R_{\mathrm{T}}(0))$ and $r(0)=1$, $r_{\mathrm{R}}'(0+) < 0$, we thus see that $\criticalpoint{\eta} = \bigO{1/|\gamma|}$, $\gamma \to -\infty$, and so in the threshold equation, see \refEquation{eqn:Alternative_form_of_the_threshold_equation}, we only need to consider $\eta$'s of $\bigO{1/|\gamma|}$. In \refEquation{eqn:Alternative_form_of_the_threshold_equation}, we have $\int_0^x \exp{(-\gamma u)} \d{u} = (1 - \exp{(-\gamma x)} ) / \gamma$. Using that $1/(\gamma B) = 1 + \bigO{1/\gamma^2}$, see \refEquation{eqn:Expansion_for_RT0_as_gamma_to_minus_infty}, we get for the right-hand side of \refEquation{eqn:Alternative_form_of_the_threshold_equation},
\begin{equation}
- \frac{1}{B} \int_0^\eta r_{\mathrm{R}}'(x) \int_0^x \e{-\gamma u} \d{u} \d{x} 
= \int_0^\eta r_{\mathrm{R}}'(x) ( 1 - \e{-\gamma x} \d{x} ) \d{x} \Bigl( 1 + \bigObig{ \frac{1}{\gamma^2} } \Bigr).
\end{equation}
Next,
\begin{equation}
r_{\mathrm{R}}'(x) 
= r_{\mathrm{R}}'(0+) + \bigObig{ \frac{1}{\gamma} }, 
\quad
1 - \e{-\gamma x} 
= \bigO{1},
\, 0 \leq x \leq \eta = \bigObig{ \frac{1}{\gamma} },
\end{equation}
and so
\begin{equation}
- \frac{1}{B} \int_0^\eta r_{\mathrm{R}}'(x) \int_0^x \e{-\gamma u} \d{u} \d{x} 
= - r_{\mathrm{R}}'(0+) \frac{ 1 - \e{-\gamma \eta} - \gamma \eta}{\gamma} + \bigObig{\frac{1}{\gamma^2}}.
\label{eqn:Expansion_of_RHS_when_gamma_to_minus_infty}
\end{equation}
Furthermore, for the left-hand side of \refEquation{eqn:Alternative_form_of_the_threshold_equation}, we have
\begin{equation}
r_{\mathrm{R}}(\eta) - R_{\mathrm{T}}(0) 
= 1 + r_{\mathrm{R}}'(0+) \eta + \bigObig{\frac{1}{\gamma^2}} - \Bigl( 1 + \frac{r_{\mathrm{L}}'(0-)}{\gamma} + \bigObig{\frac{1}{\gamma^2}} \Bigr)
= r_{\mathrm{R}}'(0+) \eta - \frac{r_{\mathrm{L}}'(0-)}{\gamma} + \bigObig{\frac{1}{\gamma^2}}.
\label{eqn:Expansion_of_LHS_when_gamma_to_minus_infty}
\end{equation}
Equating \refEquation{eqn:Expansion_of_RHS_when_gamma_to_minus_infty} and \refEquation{eqn:Expansion_of_LHS_when_gamma_to_minus_infty} and simplifying, we find
\begin{equation}
r_{\mathrm{R}}'(0+) ( 1 - \e{-\gamma \eta} ) 
= r_{\mathrm{L}}'(0-) + \bigObig{\frac{1}{\gamma}},
\end{equation}
and this gives \refEquation{eqn:Extreme_gamma__Threshold_equation_when_gamma_is_near_minus_infinity}.
}

If $r_{\mathrm{L}}(x)$ is slowly varying, the optimal threshold is approximately given by
\begin{equation}
\criticalpoint{\eta} \approx r_{\mathrm{R}}^{\gets}( r_{\mathrm{L}}(-\gamma) )
\label{eqn:Extreme_gamma__Threshold_equation_when_gamma_is_near_plus_infinity}
\end{equation}
as $\gamma \rightarrow \infty$. To see this, note that as $\gamma \rightarrow \infty$,
\begin{equation}
R_{\mathrm{T}}(0)
= \frac{ \int_{-\infty}^0 r_{\mathrm{L}}(x) \e{ -\frac{1}{2} x^2 - \gamma x } \d{x} }{ \int_{-\infty}^0 \e{ -\frac{1}{2} x^2 - \gamma x } \d{x} }
= \frac{ \e{\frac{1}{2}\gamma^2} \int_0^\infty r_{\mathrm{L}}(-x) \e{ -\frac{1}{2} ( x - \gamma )^2 } \d{x} }{ \e{\frac{1}{2}\gamma^2} \int_{0}^\infty \e{ -\frac{1}{2} ( x - \gamma )^2 } \d{x} }
\approx r_{\mathrm{L}}(-\gamma).
\end{equation}
A full analysis goes be beyond the scope of this paper, and would overly complicating our exposition. Instead, consider as an example $r_{\mathrm{L}}(x) = \exp{(bx)}$ with $b > 0$ small. We have as $\gamma \to \infty$ with exponentially small error
\begin{equation}
R_{\mathrm{T}}(0)
= \frac{\Phi(\gamma-b)}{\phi(\gamma-b)} \cdot \frac{\phi(\gamma)}{\Phi(\gamma)}
%= \frac{g(\gamma-b)}{g(\gamma)}
= \e{- b \gamma} \e{ \frac{b^2}{2} } \Bigl( 1 - b \frac{\phi(\gamma)}{\Phi(\gamma)} + \bigO{b^2} \Bigr)
= r_{\mathrm{L}}(-\gamma) ( 1 + \bigO{b^2} ).
%\approx r_{\mathrm{L}}(-\gamma) \e{ \frac{b^2}{2} }.
\label{eqn:Example_of_how_slow_varying_would_work}
\end{equation}
When for instance $r_{\mathrm{R}}(x) = \exp{(dx)}$ with $d > 0$, we get that $\criticalpoint{\eta} \approx r_{\mathrm{R}}^{\gets}( \exp{(- b \gamma + b^2 / 2 )} ) = (b/d) \gamma - b^2 / (2d)$ with exponentially small error, as $\gamma \to \infty$.
Furthermore, the right-hand side in \refEquation{eqn:Alternative_form_of_the_threshold_equation} is exponentially small as $\gamma \to \infty$, so that in good approximation the solution to the threshold equation is indeed given by \refEquation{eqn:Extreme_gamma__Threshold_equation_when_gamma_is_near_plus_infinity}.

%\begin{remark}
%\label{rem:Slow_variation_argument}
%We can be more precise about the phrase ``slowly varying'' in \refProposition{prop:Asymptotic_solutions}, at the expense of overly complicating our exposition. Instead, consider as an $r_{\mathrm{L}}(x) = \exp{(bx)}$ with $b > 0$ small. We have as $\gamma \to \infty$ with exponentially small error
%\begin{equation}
%R_{\mathrm{T}}(0)
%= \frac{\Phi(\gamma-b)}{\phi(\gamma-b)} \cdot \frac{\phi(\gamma)}{\Phi(\gamma)}
%%= \frac{g(\gamma-b)}{g(\gamma)}
%= \e{- b \gamma} \e{ \frac{b^2}{2} } \Bigl( 1 - b \frac{\phi(\gamma)}{\Phi(\gamma)} + \bigO{b^2} \Bigr)
%= r_{\mathrm{L}}(-\gamma) ( 1 + \bigO{b^2} ).
%%\approx r_{\mathrm{L}}(-\gamma) \e{ \frac{b^2}{2} }.
%\end{equation}
%When for instance $r_{\mathrm{R}}(x) = \exp{(dx)}$ with $d > 0$, we get that $\criticalpoint{\eta} \approx r_{\mathrm{R}}^{\gets}( \exp{(- b \gamma + b^2 / 2 )} ) = (b/d) \gamma - b^2 / (2d)$ with exponentially small error, as $\gamma \to \infty$.
%\end{remark}

\subsection{Explicit results for two special cases} 
\label{sec:Analyses__Exact_characterizations}

We now study the two special cases of linear and exponential revenue structures. For these cases we are able to find precise results for the $\criticalpoint{\eta}$. We demonstrate these results for some example systems, and also include the bounds and asymptotic results obtained in \refSection{sec:Analyses__General_bounds} and \refSection{sec:Analyses__Asymptotic_solutions}, respectively.

\subsubsection{Linear revenue}
\label{sec:Exact_characterizations__Linear_revenue}

We first present an explicit expression for the optimal threshold for the case of a linear revenue function,
\begin{equation}
r_{\mathrm{R}}(x) 
= \Bigl( 1 - \frac{x}{d} \Bigr) \indicator{ 0 \leq x \leq d },
\quad x \geq 0,
\end{equation}
and arbitrary $r_{\mathrm{L}}(x)$. We distinguish between $\gamma \neq 0$ and $\gamma = 0$ in \refProposition{prop:Optimal_threshold_for_linear_revenue} and \refProposition{prop:Asymptotics_of_LambertW_solutions_as_gamma_to_zero} below.

\begin{proposition}
\label{prop:Optimal_threshold_for_linear_revenue}
Assume $\gamma \neq 0$. Then
\begin{equation}
\criticalpoint{\eta} = r_0 + \frac{1}{\gamma} W\Bigl( \frac{ \gamma \e{-\gamma r_0} }{a_0} \Bigr),
\label{eqn:Optimal_threshold_for_linear_revenue}
\end{equation}
where $W$ denotes Lambert's W function, and
\begin{equation}
a_0 
= - \gamma^2 \Bigl( B + \frac{1}{\gamma} \Bigr) ,
\quad
r_0 
= \frac{ d ( B - A ) + \frac{1}{\gamma^2} }{ B + \frac{1}{\gamma} }.
\label{eqn:Definitions_of_a0_and_r0}
\end{equation}
\end{proposition}

\myProof{
It follows from \cite[Section~4]{avram_loss_2013} that $B + 1/\gamma \neq 0$ when $\gamma \neq 0$ so that $a_0$, $r_0$ in \refEquation{eqn:Definitions_of_a0_and_r0} are well-defined with $a_0 \neq 0$. From the threshold equation in \refEquation{eqn:Alternative_form_of_the_threshold_equation}, and $r_{\mathrm{R}}(\eta) = 1 - \eta / d$ when $0 \leq \eta \leq d$, we see that $\eta = \criticalpoint{\eta}$ satisfies
\begin{equation}
1 - \frac{\eta}{d} - \frac{A}{B}
= \frac{1}{d} \int_0^\eta \int_0^x \e{-\gamma u} \d{u} \d{x}.
\end{equation}
Now 
\begin{equation}
\int_0^\eta \int_0^x \e{-\gamma u} \d{u} \d{x}
= \frac{1}{\gamma} \Bigl( \eta - \frac{1-\e{-\gamma \eta}}{\gamma} \Bigr),
\end{equation}
and this yields for $\eta = \criticalpoint{\eta}$ the equation
\begin{equation}
\gamma ( \eta - r_0 ) \e{\gamma (\eta - r_0)} = \frac{\gamma}{a_0} \e{-\gamma r_0}
\label{eqn:Threshold_equation_in_Lambert_W_form_for_linear_revenue}
\end{equation}
with $a_0$ and $r_0$ given in \refEquation{eqn:Definitions_of_a0_and_r0}. The equation \refEquation{eqn:Threshold_equation_in_Lambert_W_form_for_linear_revenue} is of the form $W(z) \exp{(W(z))} = z$, which is the defining equation for Lambert's W function, and this yields the result.
}

\refProposition{prop:Optimal_threshold_for_linear_revenue} provides a connection with the developments in \cite{borgs_optimal_2014}. Furthermore, the optimal threshold $\criticalpoint{\eta}$ is readily computed from it, taking care that the branch choice for $W$ is such that the resulting $\criticalpoint{\eta}$ is positive, continuous, and increasing as a function of $\gamma$. For this matter, the following result is relevant.

\begin{proposition}
\label{prop:Asymptotics_of_LambertW_solutions_as_gamma_to_zero}
For $r_{\mathrm{R}}(x) = 1 - x/d$ with $d>0$, and arbitrary $r_{\mathrm{L}}(x)$, as $\gamma \to 0$,
\begin{equation}
\criticalpoint{\eta}
= \sqrt{ \frac{\pi}{2} + 2 d \Bigl( \sqrt{\frac{\pi}{2}} - \int_{-\infty}^0 r_{\mathrm{L}}(x) \e{-\frac{1}{2}x} \d{x} \Bigr) } - \sqrt{\frac{\pi}{2}} + \bigO{\gamma}.
\label{eqn:Asymptotics_of_LambertW_solutions_as_gamma_to_zero}
\end{equation}
\end{proposition}

\myProof{
In the threshold equation in \refEquation{eqn:Alternative_form_of_the_threshold_equation}, we set $\varepsilon = 1 - A/B$ and use $r_{\mathrm{R}}(x) = 1 - x/d$, $r_{\mathrm{R}}'(x) = - 1 / d$, to arrive at
\begin{equation}
d \varepsilon - \eta 
= \frac{1}{B} \int_0^\eta \int_0^x \e{-\gamma u} \d{u} \d{x}.
\end{equation}
Since
\begin{equation}
\int_0^\eta \int_0^x \e{-\gamma u} \d{u} \d{x}
= \int_0^\eta ( x + \bigO{ \gamma x^2 } ) \d{x} 
= \frac{1}{2} \eta^2 + \bigO{ \gamma \eta^3 },
\end{equation}
we obtain the equation
\begin{equation}
d \varepsilon - \eta = \frac{1}{2B} \eta^2 + \bigO{\gamma \eta^3}.
\label{eqn:Behavior_of_threshold_equation_for_linear_revenue_as_gamma_to_zero}
\end{equation}
Using that $\criticalpoint{\eta} < r^{\gets}( A / B) = \bigO{1}$ as $\gamma \to 0$, we find from \refEquation{eqn:Behavior_of_threshold_equation_for_linear_revenue_as_gamma_to_zero} that as $\gamma \to 0$,
\begin{equation}
\criticalpoint{\eta} 
= \sqrt{ B^2 + 2B d \varepsilon } - B + \bigO{\gamma} 
= \sqrt{ B^2 + 2d (B-A) } - B + \bigO{\gamma}.
\end{equation}
Finally \refEquation{eqn:Asymptotics_of_LambertW_solutions_as_gamma_to_zero} follows from the expansions
\begin{equation}
B = \sqrt{\frac{\pi}{2}} + \bigO{\gamma},
\quad A = \int_{-\infty}^0 r_{\mathrm{L}}(x) \e{-\frac{1}{2}x^2} \d{x} + \bigO{\gamma},
\end{equation}
as $\gamma \to 0$.
}

We may also study the regime $\gamma \rightarrow \infty$. Note that the following result coincides with the asymptotic behavior of $\eta^{\min}$ and $\eta^{\max}$ in \refProposition{prop:Threshold_control_bounds}. 
%The proof can be found in \refAppendixSection{sec:Proof_of_Asymptotics_of_LambertW_solutions_as_gamma_to_infinity}.

\begin{proposition}
\label{prop:Asymptotics_of_LambertW_solutions_as_gamma_to_infinity}
For $r_{\mathrm{L}}(x) = \e{bx}$, and $r_{\mathrm{R}}(x) = (d-x)/d$, as $\gamma \rightarrow \infty$,
\begin{equation}
\criticalpoint{\eta}
= d \Bigl( 1 - \frac{\Phi(\gamma-b) \phi(\gamma)}{\phi(\gamma-b) \Phi(\gamma)} \Bigr) + \bigObig{ \frac{1}{\gamma} \e{- \frac{1}{2} \gamma^2} }.
\end{equation}
\end{proposition}

\myProof{
The revenue structure implies that
\begin{equation}
A = \frac{\Phi(\gamma-b)}{\phi(\gamma-b)},
\quad
B = \frac{\Phi(\gamma)}{\phi(\gamma)},
\quad
\int_0^x \e{-\gamma u} \d{u} = \frac{1-\e{-\gamma x}}{\gamma}.
\end{equation}
Therefore, as $\gamma \to \infty$,
\begin{equation}
\frac{A}{B} = \frac{\Phi(\gamma-b) \phi(\gamma)}{\phi(\gamma-b) \Phi(\gamma)},
\quad
\frac{1}{B} = \bigO{ \e{-\frac{1}{2}\gamma^2} },
\quad
\int_0^x \e{-\gamma u} \d{u} = \bigObig{ \frac{1}{\gamma} }.
\end{equation}
Substituting in the threshold equation \refEquation{eqn:Alternative_form_of_the_threshold_equation}, we find that as $\gamma \to \infty$,
\begin{equation}
1 - \frac{\eta}{d} 
= \frac{\Phi(\gamma-b) \phi(\gamma)}{\phi(\gamma-b) \Phi(\gamma)} + \bigObig{ \frac{1}{\gamma} \e{ - \frac{1}{2} \gamma^2 } },
\end{equation}
which completes the proof.
}

\refFigure{fig:Figure__Threshold_as_function_of_gamma__Linear_revenue} displays $\criticalpoint{\eta}$ given in \refProposition{prop:Optimal_threshold_for_linear_revenue} as a function of $\gamma$, together with the bounds given by \refProposition{prop:Threshold_control_bounds}, 
\begin{equation}
\eta^{\max} = d \Bigl( 1 - \frac{\Phi(\gamma-b) \phi(\gamma)}{\phi(\gamma-b) \Phi(\gamma)} \Bigr),
\quad
\eta^{\min} = d \Bigl( 1 - \frac{ \Phi(\gamma-b) / \phi(\gamma-b) + \int_0^{\eta^{\max}} \e{-\gamma x} \d{x} }{ \Phi(\gamma) / \phi(\gamma) + \int_0^{\eta^{\max}} \e{-\gamma x} \d{x} } \Bigr),
\end{equation}
and asymptotic solutions of \refProposition{prop:Asymptotic_solutions}, 
\begin{equation}
\eta^{\gamma \to -\infty} = - \frac{1}{\gamma} \ln{( 1 + bd )},
\quad
\eta^{\gamma \to \infty} = d( 1 - \e{-b \gamma} ).
\end{equation}
\refFigure{fig:Figure__Threshold_as_function_of_gamma__Linear_revenue} also confirms the monotonicity of $\criticalpoint{\eta}$ in $\gamma$ established in \refProposition{prop:Optimal_threshold_is_increasing_in_gamma}. Note also the different regimes in which our approximations are valid, and that the bounds of  \refProposition{prop:Threshold_control_bounds} are tight as $\gamma \rightarrow \pm \infty$.

\begin{figure}[!hbtp]
\begin{center}
%\subfigure{
%\includegraphics[width=0.9\columnwidth]{Figure__Threshold_as_function_of_gamma__Linear_revenue}
%}
\small
%\documentclass{standalone}
%
%\usepackage{amsmath}
%\usepackage{amsthm}
%\usepackage{amssymb}
%\usepackage{siunitx}
%
%\usepackage{tikz}
%\usepackage{pgfplots}
%\usepgfplotslibrary{groupplots}
%
%\newcommand{\criticalpoint}[1]{  #1^{\textnormal{opt}} }
%
%\begin{document}
\begin{tikzpicture}
\begin{axis}[
	width=0.9\columnwidth, height=0.618*0.9\columnwidth,
    xmin=-5, xmax=5,
    ymin=0, ymax=1,
    xtick={-5,-4,...,5},
    every axis x label/.style=
        {at={(ticklabel cs:0.5)},anchor=north},
    ytick={0,1},
    every axis y label/.style=
        {at={(ticklabel cs:0.5)},rotate=90,anchor=south},
    scaled ticks=true,
	xlabel={$\gamma$},
	]

%%%%%%%%%%%%%%%%%%%%%%%%%%%
%% The threshold graphs. %%
%%%%%%%%%%%%%%%%%%%%%%%%%%%

% Optimal threshold, etaOpt.
\addplot[color=black, thick, mark=none] plot coordinates {
(-5.,0.115494) (-4.89899,0.11738) (-4.79798,0.119326) (-4.69697,0.121334) (-4.59596,0.123407) (-4.49495,0.125548) (-4.39394,0.127762) (-4.29293,0.13005) (-4.19192,0.132417) (-4.09091,0.134867) (-3.9899,0.137403) (-3.88889,0.140031) (-3.78788,0.142754) (-3.68687,0.145578) (-3.58586,0.148509) (-3.48485,0.151551) (-3.38384,0.15471) (-3.28283,0.157994) (-3.18182,0.161409) (-3.08081,0.164962) (-2.9798,0.168661) (-2.87879,0.172514) (-2.77778,0.176531) (-2.67677,0.18072) (-2.57576,0.185093) (-2.47475,0.189659) (-2.37374,0.19443) (-2.27273,0.19942) (-2.17172,0.20464) (-2.07071,0.210106) (-1.9697,0.215832) (-1.86869,0.221835) (-1.76768,0.228132) (-1.66667,0.234741) (-1.56566,0.241682) (-1.46465,0.248975) (-1.36364,0.256642) (-1.26263,0.264707) (-1.16162,0.273194) (-1.06061,0.282128) (-0.959596,0.291536) (-0.858586,0.301445) (-0.757576,0.311884) (-0.656566,0.32288) (-0.555556,0.334463) (-0.454545,0.346661) (-0.353535,0.359501) (-0.252525,0.373009) (-0.151515,0.387206) (-0.0505051,0.402112) (0.0505051,0.41774) (0.151515,0.434098) (0.252525,0.451184) (0.353535,0.468987) (0.454545,0.487484) (0.555556,0.50664) (0.656566,0.526403) (0.757576,0.546708) (0.858586,0.567471) (0.959596,0.588596) (1.06061,0.609968) (1.16162,0.631463) (1.26263,0.652946) (1.36364,0.674278) (1.46465,0.695319) (1.56566,0.715931) (1.66667,0.735987) (1.76768,0.755371) (1.86869,0.773983) (1.9697,0.791742) (2.07071,0.808584) (2.17172,0.824467) (2.27273,0.839367) (2.37374,0.853277) (2.47475,0.866204) (2.57576,0.878171) (2.67677,0.889207) (2.77778,0.899355) (2.87879,0.908658) (2.9798,0.917166) (3.08081,0.924931) (3.18182,0.932005) (3.28283,0.93844) (3.38384,0.944285) (3.48485,0.94959) (3.58586,0.954401) (3.68687,0.958759) (3.78788,0.962706) (3.88889,0.966279) (3.9899,0.969512) (4.09091,0.972437) (4.19192,0.975082) (4.29293,0.977475) (4.39394,0.979638) (4.49495,0.981593) (4.59596,0.983361) (4.69697,0.984959) (4.79798,0.986404) (4.89899,0.98771) (5.,0.988891) 
};

% Upper bound, etaMax.
\addplot[color=black, thick, style=dashed, mark=none] plot coordinates {
(-5.,0.157828) (-4.89899,0.160246) (-4.79798,0.162735) (-4.69697,0.165298) (-4.59596,0.167938) (-4.49495,0.170658) (-4.39394,0.173462) (-4.29293,0.176354) (-4.19192,0.179336) (-4.09091,0.182414) (-3.9899,0.185591) (-3.88889,0.188872) (-3.78788,0.192262) (-3.68687,0.195765) (-3.58586,0.199388) (-3.48485,0.203135) (-3.38384,0.207012) (-3.28283,0.211026) (-3.18182,0.215182) (-3.08081,0.219489) (-2.9798,0.223953) (-2.87879,0.228581) (-2.77778,0.233382) (-2.67677,0.238365) (-2.57576,0.243537) (-2.47475,0.24891) (-2.37374,0.254491) (-2.27273,0.260293) (-2.17172,0.266325) (-2.07071,0.272599) (-1.9697,0.279127) (-1.86869,0.285922) (-1.76768,0.292996) (-1.66667,0.300363) (-1.56566,0.308036) (-1.46465,0.31603) (-1.36364,0.324359) (-1.26263,0.333039) (-1.16162,0.342084) (-1.06061,0.35151) (-0.959596,0.361332) (-0.858586,0.371565) (-0.757576,0.382222) (-0.656566,0.393319) (-0.555556,0.404866) (-0.454545,0.416877) (-0.353535,0.429359) (-0.252525,0.442321) (-0.151515,0.455766) (-0.0505051,0.469697) (0.0505051,0.48411) (0.151515,0.498998) (0.252525,0.514349) (0.353535,0.530146) (0.454545,0.546363) (0.555556,0.56297) (0.656566,0.579929) (0.757576,0.597194) (0.858586,0.614712) (0.959596,0.632423) (1.06061,0.65026) (1.16162,0.668149) (1.26263,0.686012) (1.36364,0.703766) (1.46465,0.721329) (1.56566,0.738615) (1.66667,0.755541) (1.76768,0.772029) (1.86869,0.788005) (1.9697,0.803403) (2.07071,0.818165) (2.17172,0.832244) (2.27273,0.845604) (2.37374,0.858218) (2.47475,0.870073) (2.57576,0.881164) (2.67677,0.891496) (2.77778,0.901084) (2.87879,0.90995) (2.9798,0.918121) (3.08081,0.925628) (3.18182,0.932509) (3.28283,0.938799) (3.38384,0.94454) (3.48485,0.949768) (3.58586,0.954523) (3.68687,0.958843) (3.78788,0.962763) (3.88889,0.966317) (3.9899,0.969538) (4.09091,0.972454) (4.19192,0.975093) (4.29293,0.977481) (4.39394,0.979642) (4.49495,0.981596) (4.59596,0.983363) (4.69697,0.984961) (4.79798,0.986405) (4.89899,0.987711) (5.,0.988891) 
};

% Lower bound, etaMin.
\addplot[color=black, thick, style=dashed, mark=none] plot coordinates {
(-5.,0.0702607) (-4.89899,0.0715793) (-4.79798,0.0729461) (-4.69697,0.0743634) (-4.59596,0.0758341) (-4.49495,0.0773611) (-4.39394,0.0789474) (-4.29293,0.0805964) (-4.19192,0.0823118) (-4.09091,0.0840972) (-3.9899,0.0859568) (-3.88889,0.0878951) (-3.78788,0.0899168) (-3.68687,0.0920271) (-3.58586,0.0942313) (-3.48485,0.0965355) (-3.38384,0.0989461) (-3.28283,0.10147) (-3.18182,0.104115) (-3.08081,0.106888) (-2.9798,0.109799) (-2.87879,0.112858) (-2.77778,0.116073) (-2.67677,0.119457) (-2.57576,0.123022) (-2.47475,0.126781) (-2.37374,0.130747) (-2.27273,0.134937) (-2.17172,0.139367) (-2.07071,0.144055) (-1.9697,0.149022) (-1.86869,0.154288) (-1.76768,0.159877) (-1.66667,0.165814) (-1.56566,0.172126) (-1.46465,0.178843) (-1.36364,0.185997) (-1.26263,0.19362) (-1.16162,0.201751) (-1.06061,0.210426) (-0.959596,0.219688) (-0.858586,0.229578) (-0.757576,0.240142) (-0.656566,0.251427) (-0.555556,0.263477) (-0.454545,0.276341) (-0.353535,0.290063) (-0.252525,0.304686) (-0.151515,0.320249) (-0.0505051,0.336783) (0.0505051,0.354312) (0.151515,0.372848) (0.252525,0.392389) (0.353535,0.412916) (0.454545,0.434391) (0.555556,0.456752) (0.656566,0.479917) (0.757576,0.503776) (0.858586,0.528199) (0.959596,0.55303) (1.06061,0.5781) (1.16162,0.603223) (1.26263,0.628208) (1.36364,0.652864) (1.46465,0.677005) (1.56566,0.700462) (1.66667,0.723083) (1.76768,0.744741) (1.86869,0.765337) (1.9697,0.784796) (2.07071,0.803074) (2.17172,0.82015) (2.27273,0.836026) (2.37374,0.850721) (2.47475,0.864273) (2.57576,0.876727) (2.67677,0.888142) (2.77778,0.898576) (2.87879,0.908095) (2.9798,0.916764) (3.08081,0.924647) (3.18182,0.931806) (3.28283,0.938302) (3.38384,0.944191) (3.48485,0.949526) (3.58586,0.954358) (3.68687,0.958731) (3.78788,0.962688) (3.88889,0.966267) (3.9899,0.969504) (4.09091,0.972432) (4.19192,0.975079) (4.29293,0.977473) (4.39394,0.979636) (4.49495,0.981592) (4.59596,0.983361) (4.69697,0.984959) (4.79798,0.986404) (4.89899,0.98771) (5.,0.988891) 
};

% Asymptotic regime gamma to infinity, etaGammaToInfinity.
\addplot[color=black, thick, style=dotted, mark=none] plot coordinates {
(0.0505051,0.0492509) (0.151515,0.140595) (0.252525,0.223163) (0.353535,0.297799) (0.454545,0.365264) (0.555556,0.426247) (0.656566,0.481371) (0.757576,0.531198) (0.858586,0.576239) (0.959596,0.616952) (1.06061,0.653754) (1.16162,0.68702) (1.26263,0.71709) (1.36364,0.744271) (1.46465,0.76884) (1.56566,0.791049) (1.66667,0.811124) (1.76768,0.829271) (1.86869,0.845674) (1.9697,0.860501) (2.07071,0.873903) (2.17172,0.886018) (2.27273,0.896969) (2.37374,0.906868) (2.47475,0.915816) (2.57576,0.923904) (2.67677,0.931215) (2.77778,0.937823) (2.87879,0.943797) (2.9798,0.949197) (3.08081,0.954078) (3.18182,0.95849) (3.28283,0.962478) (3.38384,0.966083) (3.48485,0.969342) (3.58586,0.972287) (3.68687,0.97495) (3.78788,0.977356) (3.88889,0.979532) (3.9899,0.981498) (4.09091,0.983276) (4.19192,0.984883) (4.29293,0.986335) (4.39394,0.987648) (4.49495,0.988835) (4.59596,0.989907) (4.69697,0.990877) (4.79798,0.991754) (4.89899,0.992546) (5.,0.993262) 
};

% Asymptotic regime gamma to minus infinity, etaGammaToMinusInfinity.
\addplot[color=black, thick, style=dotted, mark=none] plot coordinates {
(-5.,0.138629) (-4.89899,0.141488) (-4.79798,0.144466) (-4.69697,0.147573) (-4.59596,0.150817) (-4.49495,0.154206) (-4.39394,0.157751) (-4.29293,0.161463) (-4.19192,0.165353) (-4.09091,0.169436) (-3.9899,0.173725) (-3.88889,0.178238) (-3.78788,0.182991) (-3.68687,0.188004) (-3.58586,0.1933) (-3.48485,0.198903) (-3.38384,0.204841) (-3.28283,0.211143) (-3.18182,0.217846) (-3.08081,0.224989) (-2.9798,0.232615) (-2.87879,0.240777) (-2.77778,0.249533) (-2.67677,0.258949) (-2.57576,0.269104) (-2.47475,0.280088) (-2.37374,0.292007) (-2.27273,0.304985) (-2.17172,0.31917) (-2.07071,0.334739) (-1.9697,0.351905) (-1.86869,0.370927) (-1.76768,0.392123) (-1.66667,0.415888) (-1.56566,0.44272) (-1.46465,0.473252) (-1.36364,0.508308) (-1.26263,0.548973) (-1.16162,0.596709) (-1.06061,0.653539) (-0.959596,0.722332) (-0.858586,0.807313) (-0.757576,0.914954) (-0.656566,1.05572) (-0.555556,1.24766) (-0.454545,1.52492) (-0.353535,1.96062) (-0.252525,2.74486) (-0.151515,4.57477) (-0.0505051,13.7243) 
};

\end{axis}

%\draw (4.5,1.5) node [fill=none, draw=none, text=black, opacity=1, text opacity=1] {$\criticalpoint{\eta}$};
\draw (4.75,1.05) node [fill=none, draw=none, text=black, opacity=1, text opacity=1] {$\eta^{\min}$};
\draw (5.25,3) node [fill=none, draw=none, text=black, opacity=1, text opacity=1] {$\eta^{\max}$};
\draw (7.75,5.75) node [fill=none, draw=none, text=black, opacity=1, text opacity=1] {$\eta^{\gamma \to \infty}$};
\draw (4,4) node [fill=none, draw=none, text=black, opacity=1, text opacity=1] {$\eta^{\gamma \to -\infty}$};

\end{tikzpicture}
%\end{document}
\end{center}
\vspace{-1em}
\caption{\textrm{The optimal threshold $\criticalpoint{\eta}$, its bounds $\eta^{\min}$, $\eta^{\max}$, and its approximations $\eta^{\gamma \to \pm \infty}$, all as a function of $\gamma$, when $r_{\mathrm{L}}(x) = \exp{(bx)}$, $r_{\mathrm{R}}(x) = (d-x)/d$, and $b = d = 1$. The curve for the optimal threshold has been produced with \refEquation{eqn:Optimal_threshold_for_linear_revenue}.}}
\label{fig:Figure__Threshold_as_function_of_gamma__Linear_revenue}
\end{figure}

\subsubsection{Exponential revenue} 
\label{sec:Exact_characterizations__Exponential_revenue}

Consider $r_{\mathrm{L}}(x)$ arbitrary, and let $r_{\mathrm{R}}(x) = \exp{(-\delta x)}$ for $x \geq 0$, with $\delta > 0$. First, we will consider what happens asymptotically as $\delta \downarrow 0$ in the case $\gamma = 0$, which should be comparable to the case in \refProposition{prop:Asymptotics_of_LambertW_solutions_as_gamma_to_zero}. Then, we consider the case $\gamma = -\delta$, which like the linear revenue structure has a Lambert W solution. Finally, we consider what happens asymptotically when $\varepsilon = 1 - R_{\mathrm{T}}(x) > 0$ is small, and we check our results in the specific cases $\gamma = -2\delta$, $-\delta/2$ and $\gamma = \delta$, which have explicit solutions.

\begin{proposition}
\label{prop:Exponential_revenue__Asymptotics_as_delta_goes_to_zero}
For $\gamma = 0$, as $\delta \downarrow 0$,
\begin{equation}
\criticalpoint{\eta}
= \sqrt{ \frac{2(B-A)}{\delta} } - \frac{2A + B}{3}  + \bigObig{ \sqrt{\delta} }
= \sqrt{ \frac{2}{\delta} \Bigl( \sqrt{\frac{\pi}{2}} - \int_{-\infty}^0 r(x) \e{-\frac{1}{2} x^2} \d{x} \Bigr) } + \bigO{1}.
\label{eqn:Exponential_revenue__Solution_when_gamma_is_0_and_delta_to_zero}
\end{equation}
\end{proposition}

\myProof{
When $\gamma = 0$, the threshold equation reads
\begin{equation}
\e{\delta \eta} 
= 1 + \frac{\delta(B-A)}{1 + A \delta} + \frac{ \delta \eta }{ 1 + A \delta },
\label{eqn:Exponential_case__Threshold_equation_when_gamma_is_0}
\end{equation}
which follows from \refEquation{eqn:Threshold_equation_rewritten} with $\gamma = 0$. With $\delta > 0$, $\eta > 0$, the left-hand side of \refEquation{eqn:Exponential_case__Threshold_equation_when_gamma_is_0} exceeds $1 + \delta \eta + \delta^2 \eta^2 / 2$, while the right-hand side is exceeded by $1 + \delta \eta + \delta ( B - A )$. Therefore, the left-hand side of \refEquation{eqn:Exponential_case__Threshold_equation_when_gamma_is_0} exceeds the right-hand side if $\eta > \eta_*$, where $\eta_* = \sqrt{ 2(B-A) / \delta }$. This implies that $\criticalpoint{\eta} \leq \eta^*$, and so we restrict attention to $0 \leq \eta \leq \eta_* = \bigO{ 1/\sqrt{\delta} }$ when considering \refEquation{eqn:Exponential_case__Threshold_equation_when_gamma_is_0}. Expanding both sides of \refEquation{eqn:Exponential_case__Threshold_equation_when_gamma_is_0} gives
\begin{equation}
1 + \delta \eta + \frac{1}{2} \delta^2 \eta^2 + \frac{1}{6} \delta^3 \eta^3 + \bigO{ \delta^4 \eta^4 }
= 1 + \delta(B-A) - \delta^2 A(B-A) + \bigO{\delta^3} + \delta \eta - \delta^2 A \eta + \bigO{ \delta^3 \eta }.
\label{eqn:Expansion_of_the_threshold_equation_for_gamma_0_and_exponential_case}
\end{equation}
Cancelling the terms $1 + \delta \eta$ at both sides of \refEquation{eqn:Expansion_of_the_threshold_equation_for_gamma_0_and_exponential_case}, and dividing by $\delta^2 / 2$ while remembering that $\eta = \bigO{ 1 / \sqrt{\delta} }$, we get
\begin{equation}
\eta^2 
= \frac{2(B-A)}{\delta} - 2 \eta A - \frac{1}{3} \eta^3 \delta + \bigO{1}.
\end{equation}
Therefore, 
\begin{equation}
\eta
= \eta_* \Bigl( 1 - \frac{A \delta}{B-A} \eta - \frac{\delta^2}{6(B-A)} \eta^3 + \bigO{\delta} \Bigr)^{\frac{1}{2}}
= \eta_* ( 1 + \bigO{\sqrt{\delta}} ).
\label{eqn:eta_in_terms_of_eta_star}
\end{equation}
Thus $\eta = \eta_* + \bigO{1}$, and inserting this in the right-hand side of the middle member in \refEquation{eqn:eta_in_terms_of_eta_star} yields
\begin{align}
\eta 
&
= \eta_* \Bigl( 1 - \frac{A \delta}{B-A} \eta_* - \frac{\delta^2}{6(B-A)} \eta_*^3 + \bigO{\delta} \Bigr)^{\frac{1}{2}}
\nonumber \\ &
= \eta_* - \frac{A \delta}{2(B-A)} \eta_*^2 - \frac{\delta^2}{12(B-A)} \eta_*^4 + \bigO{\sqrt{\delta}}
\nonumber \\ &
= \eta_* - \frac{2A+B}{3} + \bigO{\sqrt{\delta}},
\end{align}
which is the result \refEquation{eqn:Exponential_revenue__Solution_when_gamma_is_0_and_delta_to_zero}.
}

\refFigure{fig:Figure__Threshold_as_function_of_Delta} draws for $\gamma = 0$ a comparison between $r_{\mathrm{R}}(x) = \exp{(-\Delta x)}$ and $r_{\mathrm{R}}(x) = 1 - \Delta x$. As expected, we see agreement when $\Delta \downarrow 0$, and for larger $\Delta$ the exponential revenue leads to slightly larger $\criticalpoint{\eta}$ compared with linear revenues.
%%%%%%%%%%%%%%%%%%
%%%% Old text. %%%
%%%%%%%%%%%%%%%%%%
%The exponential revenue and linear revenue resemble one another when $1/d$ and $\delta$ are equal and small positive. This is particularly evident when comparing \refProposition{prop:Asymptotics_of_LambertW_solutions_as_gamma_to_zero} with \refProposition{prop:Exponential_revenue__Asymptotics_as_delta_goes_to_zero}. To illustrate this point further, we plot the optimal threshold in the linear case (given explicitly in \refProposition{prop:Asymptotics_of_LambertW_solutions_as_gamma_to_zero}), together with the optimal threshold in the exponential case (determined numerically) in \refFigure{fig:Figure__Threshold_as_function_of_Delta} for $\gamma = 0$. We have also included the leading-order asymptotic behavior, as ascertained in \refProposition{prop:Exponential_revenue__Asymptotics_as_delta_goes_to_zero}.

\begin{figure}[!hbtp]
\begin{center}
%\subfigure{
%\includegraphics[width=0.76\columnwidth]{Figure__Threshold_as_function_of_Delta}
%}
\small
\input{Figure__Threshold_as_function_of_Delta}
\end{center}
\vspace{-1em}
\caption{\textrm{The optimal threshold $\criticalpoint{\eta}$ in the exponential revenue case $r_{\mathrm{R}}(x) = \exp{(-\Delta x)}$, and in the linear revenue case $r_{\mathrm{R}}(x) = 1 - \Delta x$, as $\Delta \downarrow 0$. In both cases, $r_{\mathrm{R}}'(0+) = - \Delta$. The leading-order behavior established in \refProposition{prop:Exponential_revenue__Asymptotics_as_delta_goes_to_zero} is also included.}}
\label{fig:Figure__Threshold_as_function_of_Delta}
\end{figure}

When $\gamma = -\delta$, the threshold equation becomes
\begin{equation}
\e{ - \delta \eta } - R_{\mathrm{T}}(0)
= \frac{\phi(-\delta)}{\Phi(-\delta)} \Bigl( \eta - \frac{1 - \e{-\delta \eta}}{\delta} \Bigr),
\end{equation}
or equivalently,
\begin{equation}
\e{ - \delta \eta } = \frac{1}{ \frac{\Phi(-\delta)}{\phi(-\delta)} - \frac{1}{\delta} } \Bigl( \eta - \frac{1}{\delta} + \frac{\Phi(-\delta)}{\phi(-\delta)} R_{\mathrm{T}}(0) \Bigr),
\end{equation}
and the solution may again be expressed in terms the Lambert W function.

\begin{proposition}
\label{prop:Exponential_revenue__LambertW_solution}
When $\gamma = - \delta$, $\criticalpoint{\eta} = r_0 + (1/\delta) W( \delta \e{-\delta r_0} / a_0 )$, with
\begin{equation}
a_0 = \frac{1}{ \frac{\Phi(\gamma)}{\phi(\gamma)} - \frac{1}{\delta} },
\quad 
r_0 = \frac{1}{\delta} - \frac{\Phi(\gamma)}{\phi(\gamma)} R_{\mathrm{T}}(0).
\end{equation}
\end{proposition}

\myProof{
Immediate, since the standard form is $\e{-\delta \eta} = a_0 ( \eta - r_0 )$.
}

In case $\alpha = (\gamma+\delta) / \delta \neq 0, 1$, the threshold equation is given by, see \refEquation{eqn:Alternative_form_of_the_threshold_equation},
\begin{equation}
\e{ - \delta \eta } - R_{\mathrm{T}}(0) 
= \frac{\delta}{B} \int_0^\eta \e{-\delta x} \frac{1 - \e{-\gamma x}}{\gamma} \d{x}
= \frac{\delta}{B\gamma} \Bigl( \frac{1 - \e{-\delta \eta}}{\delta} - \frac{1 - \e{-(\gamma+\delta)\eta}}{\gamma+\delta} \Bigr),
\label{eqn:Threshold_equation_for_exponential_revenue}
\end{equation}
After setting $z = \e{-\delta \eta} \in (0,1]$, \refEquation{eqn:Threshold_equation_for_exponential_revenue} takes the form
\begin{equation}
z - R_{\mathrm{T}}(0) = \frac{1}{\gamma B} \Bigl( 1 - z - \frac{1 - z^\alpha}{\alpha} \Bigr).
\label{eqn:Threshold_equation_for_exponential_revenue_in_terms_of_z}
\end{equation}

Observe that the factor $1 / \gamma B$ is positive when $\alpha > 1$, and negative when $\alpha < 1$. For values $\alpha = -1$, $1/2$, and $2$, an explicit solution can be found in terms of the square-root function, see \refProposition{prop:Exponential_revenue__Square_root_solutions} in \refAppendixSection{sec:Appendix__Explicit_solutions_for_exponential_revenue}. In all other cases, the solution is more involved. In certain regimes, however, a solution in terms of an infinite power series can be obtained, see \refProposition{prop:Explicit_characterization_in_exponential_case_when_R0_is_approximately_one} in \refAppendixSection{sec:Appendix__Explicit_solutions_for_exponential_revenue}.

For illustrative purposes, we again plot the optimal threshold $\criticalpoint{\eta}$ as a function of $\gamma$. It has been determined by numerically solving the threshold equation, and is plotted together with the bounds given by \refProposition{prop:Threshold_control_bounds}, 
\begin{equation}
\eta^{\max} = - \frac{1}{d} \ln{ \Bigl( \frac{ g(\gamma-b) }{ g(\gamma) } \Bigr) },
\quad
\eta^{\min} = - \frac{1}{d} \ln{ \Bigl( \frac{ g(\gamma-b) + \int_0^{\eta^{\max}} \e{-\gamma x} \d{x} }{ g(\gamma) + \int_0^{\eta^{\max}} \e{-\gamma x} \d{x} } \Bigr) },
\end{equation}
and asymptotic solutions of \refProposition{prop:Asymptotic_solutions}, 
\begin{equation}
\eta^{\gamma \to -\infty} = - \frac{1}{\gamma} \ln{ \Bigl( 1 + \frac{b}{d} \Bigr) },
\quad
\eta^{\gamma \to \infty} = \frac{b \gamma}{d}, 
\end{equation}
in \refFigure{fig:Figure__Threshold_as_function_of_gamma__Exponential_revenue}. Similar to \refFigure{fig:Figure__Threshold_as_function_of_gamma__Linear_revenue}, \refFigure{fig:Figure__Threshold_as_function_of_gamma__Exponential_revenue} also illustrates the monotonicity of $\criticalpoint{\eta}$ in $\gamma \in \realNumbers$, the different regimes our approximations and bounds are valid, and how our bounds are tight as $\gamma \rightarrow \infty$. We have also indicated the analytical solutions for $\alpha = -1$, $0$, $1/2$, and $2$, as provided by \refProposition{prop:Exponential_revenue__LambertW_solution} and \refProposition{prop:Exponential_revenue__Square_root_solutions} in \refAppendixSection{sec:Appendix__Explicit_solutions_for_exponential_revenue}. The asymptotic width $1/2$ of the gap between the graphs of $\criticalpoint{\eta}$ and $\eta^{\gamma \to \infty}$ is consistent with the refined asymptotics of $\criticalpoint{\eta}$ as given below \refEquation{eqn:Example_of_how_slow_varying_would_work}, case $b = d = 1$.

\begin{figure}[!hbtp]
\begin{center}
%\subfigure{
%\includegraphics[width=0.9\columnwidth]{Figure__Threshold_as_function_of_gamma__Exponential_revenue}
%}
\small
\input{Figure__Threshold_as_function_of_gamma__Exponential_revenue}
\end{center}
\vspace{-1em}
\caption{\textrm{The optimal threshold $\criticalpoint{\eta}$, its bounds $\eta^{\min}$, $\eta^{\max}$, and its approximations $\eta^{\gamma \to \pm \infty}$, all as a function of $\gamma$, when $r_{\mathrm{L}}(x) = \exp{(bx)}$, $r_{\mathrm{R}}(x) = \exp{(-dx)}$, and $b = d = 1$. The analytical solutions for $\alpha = -1$, $0$, $1/2$, and $2$ provided by \refProposition{prop:Exponential_revenue__LambertW_solution} and \refProposition{prop:Exponential_revenue__Square_root_solutions} are also indicated. The curve for the optimal threshold has been produced by numerically solving the threshold equation.}}
\label{fig:Figure__Threshold_as_function_of_gamma__Exponential_revenue}
\end{figure}

\section{Optimality of threshold policies}
\label{sec:Optimality_of_threshold_policies}

We now present a proof of \refProposition{prop:Threshold_control_is_optimal_form_of_control}, the cornerstone for this paper that says that threshold policies are optimal, and that the optimal threshold satisfies the threshold equation. We first present in \refSection{sec:Variational_argument} a variational argument that gives an insightful way to derive \refProposition{prop:Threshold_control_is_optimal_form_of_control} heuristically. Next, we present the formal proof of \refProposition{prop:Threshold_control_is_optimal_form_of_control} in \refSection{sec:Proof_that_threshold_controls_are_optimal} using Hilbert-space theory.

\subsection{Heuristic based on a variational argument}
\label{sec:Variational_argument}

For threshold controls $f(x) = \indicator{ 0 \leq x < \eta }$ with $\eta \in \positiveRealNumbers$, the \gls{QED} limit of the long-term revenue \refEquation{eqn:QED_limit_of_revenue_rate} becomes \refEquation{eqn:Long_term_revenue__Threshold}. The optimal threshold $\criticalpoint{\eta}$ can be found by equating
\begin{equation}
\frac{\d{R}}{\d{\eta}}
= \frac{ \bigl( B + \frac{1 - \e{-\gamma \eta} }{\gamma} \bigr) r(\eta) \e{-\gamma \eta} - \bigl( A + \int_0^\eta r(x) \e{-\gamma x} \d{x} \bigr) \e{ - \gamma \eta } }{ \bigl( B + \frac{1 - \e{-\gamma \eta} }{\gamma} \bigr)^2 } 
= \frac{ r(\eta) - R_{\mathrm{T}}(\eta) }{ \e{\gamma \eta} \Bigl( B + \frac{1 - \e{-\gamma \eta} }{\gamma} \Bigr) }
\label{eqn:Optimal_threshold_within_the_class_of_threshold_controls}
\end{equation}
to zero, which shows that the optimal threshold $\criticalpoint{\eta}$ solves the threshold equation \refEquation{eqn:Threshold_equation}, i.e.\ $r(\eta) = R_{\mathrm{T}}(\eta)$.

For any piecewise continuous function $g$ on $[0,\infty)$ that is \emph{admissible}, i.e.\ such that $0 \leq f + \varepsilon g \leq 1$ and $\int_0^\infty ( f + \varepsilon g ) \e{-\gamma x} \d{x} < \infty$ for sufficiently small $\varepsilon$, define
\begin{equation}
\delta R(f ; g) 
= \lim_{\varepsilon \downarrow 0} \frac{ R(f + \varepsilon g) - R(f) }{ \varepsilon }.
\label{eqn:Definition__Gatteaux_derivative}
\end{equation}
We call \refEquation{eqn:Definition__Gatteaux_derivative} the functional derivative of $f$ with increment $g$, which can loosely be interpreted as a derivative of $f$ in the direction of $g$, see \cite{luenberger_optimization_1969} for background. Substituting \refEquation{eqn:QED_limit_of_revenue_rate} into \refEquation{eqn:Definition__Gatteaux_derivative} yields
\begin{equation}
\delta R(f ; g)
= \frac{ \bigl( B + \int_0^\infty f \e{-\gamma x} \d{x} \bigr) \int_0^\infty r g \e{-\gamma x} \d{x} - \bigl( A + \int_0^\infty r f \e{-\gamma x} \d{x} \bigr) \int_0^\infty g \e{-\gamma x} \d{x} }{ \bigl( B  + \int_0^\infty f \e{-\gamma x} \d{x} \bigr)^2 }
\label{eqn:Gatteax_derivative_of_Rf}
\end{equation}
Rewriting \refEquation{eqn:Gatteax_derivative_of_Rf} gives
\begin{equation}
\delta R(f ; g)
= \frac{ \int_{0}^\infty g(x) \e{-\gamma x} \bigl[ r(x) - R(f) \bigr] \d{x} }{ B + \int_0^\infty f(x) \e{-\gamma x} \d{x} }.
\label{eqn:Gatteaux_differential_of_Rf__Simplified}
\end{equation} 

We can now examine the effect of small perturbations $\varepsilon g$ towards (or away from) policies $f$ by studying the sign of \refEquation{eqn:Gatteaux_differential_of_Rf__Simplified}. Specifically, it can be shown that for every perturbation $g$ applied to the optimal threshold policy of \refProposition{prop:Threshold_control_is_optimal_form_of_control}, $\delta R( \criticalpoint{f} ; g ) \leq 0$, indicating that these threshold policies are locally optimal. Moreover, it can be shown that for any other control $f$, a perturbation exists so that $\delta R( f ; g ) > 0$. Such other controls are therefore not locally optimal. Assuming the existence of an optimizer, these observations thus indeed indicate that the threshold control in \refProposition{prop:Threshold_control_is_optimal_form_of_control} is optimal. We note that these observations crucially depend on the sign of $r(x) - R(f)$, as can be seen from \refEquation{eqn:Gatteaux_differential_of_Rf__Simplified}. It is in fact the threshold equation \refEquation{eqn:Threshold_equation} that specifies the point where a sign change occurs. 

Note that while these arguments support \refProposition{prop:Threshold_control_is_optimal_form_of_control}, this section does not constitute a complete proof. In particular the existence of optimizers still needs to be established.

\subsection{Formal proof of \refProposition{prop:Threshold_control_is_optimal_form_of_control}}
\label{sec:Proof_that_threshold_controls_are_optimal}

In the formal proof of \refProposition{prop:Threshold_control_is_optimal_form_of_control} that now follows, we start by proving that there exist maximizers in \refSection{sec:Proof_of_optimal_threshold_controls__Existence_of_allowed_maximizers}. This ensures that our maximization problem is well-defined. In \refSection{sec:Proof_of_optimal_threshold_controls__Necessary_condition_for_maximizers}, we then derive necessary conditions for maximizers by perturbing the control towards (or away from) a threshold policy, as alluded to before, and in a formal manner using measure theory. Finally, we characterize in \refSection{sec:Proof_of_optimal_threshold_controls__Characterization_of_maximizers} the maximizers, by formally discarding pathological candidates.

With $r: \realNumbers \rightarrow \positiveRealNumbers$ a smooth function, nonincreasing to $0$ as $x \rightarrow \pm \infty$, and $\gamma \in \realNumbers$, recall that we are considering the maximization of the functional \refEquation{eqn:QED_limit_of_revenue_rate} with $f : \positiveRealNumbers \rightarrow [0,1]$ is measurable and with $g(x) = f(x) \e{-\gamma x} \in L^1( \positiveRealNumbers )$. We do not assume $f$ to be nonincreasing. Recall that
$
A = \int_{-\infty}^0 r(x) \exp{(- \frac{1}{2} x^2 - \gamma x )} \d{x} > 0
$,
$
B = \Phi(\gamma) / \phi(\gamma) > 0
$,
and let $b(x) = \e{-\gamma x}$ for $x \geq 0$. Then write $R(f)$ as 
\begin{equation}
R(f) 
= \frac{ A + \int_0^\infty r(x) g(x) \d{x} }{ B + \int_0^\infty g(x) \d{x} } = L(g),
\end{equation}
which is considered for all $g \in L^1( \positiveRealNumbers )$ such that $0 \leq g(x) \leq b(x)$ for $0 \leq x < \infty$. The objective is to maximize $L(g)$ over all such allowed $g$.

For notational convenience, write
\begin{equation}
L(g) 
= \frac{A}{B} \Bigl( 1 + \frac{ \int_0^\infty s(x) g(x) \d{x} }{ 1 + \int_0^\infty S g(x) \d{x} } \Bigr),
\label{eqn:Definition_of_Lg}
\end{equation}
where 
\begin{equation}
s(x) 
= \frac{r(x)}{A} - \frac{1}{B},
\quad 
S 
= \frac{1}{B}.
\end{equation}

Recall that $r(x)$ is nonincreasing, implying that $s(x) \leq s(0)$ for all $x \geq 0$. When $s(0) \leq 0$, the maximum of \refEquation{eqn:Definition_of_Lg} thus equals $A/B$, and is assumed by all allowed $g$ that vanish outside the interval $[0, \sup\{ x \in \positiveRealNumbers | s(x) = 0 \} ]$. When $s(0) > 0$, define
\begin{equation}
x_0 = \inf \{ x \in \positiveRealNumbers | s(x) = 0 \}, 
\label{eqn:Definition__Set_x0}
\end{equation}
which is positive and finite by smoothness of $s$ and $r(x) \rightarrow 0$ as $x \rightarrow \infty$. Note that the set $\{ x \in \positiveRealNumbers | s(x) = 0 \}$ consists of a single point when $r(x)$ is \emph{strictly} decreasing as long as $r(x) > 0$. But even if $r(x)$ is not strictly decreasing, we have $s(x) \leq 0$ for $x \geq x_0$. Because $g(x) \geq 0$ implies that 
\begin{equation}
\int_{x_0}^\infty s(x) g(x) \d{x} 
\leq 0 
\leq \int_{x_0}^\infty S g(x) \d{x},
\label{eqn:Bounds__On_integral_over_sx_gx}
\end{equation}
we have 
\begin{equation}
L( g \indicator{ x \in [0, x_0) } ) \geq L(g)
\end{equation} 
for all $g$. We may therefore restrict attention to allowed $g$ supported on $[0,x_0]$. Such a $g$ can be extended to any allowed function supported on $[0,\sup\{ x \in \positiveRealNumbers | s(x) = 0 \}$ without changing the value $L(g)$. Therefore, we shall instead maximize
\begin{equation}
J(g)
= \frac{ \int_0^{x_0} s(x) g(x) \d{x} }{ 1 + \int_0^{x_0} S g(x) \d{x} }
\label{eqn:Definition_of_Jg}
\end{equation}
over all $g \in L^1([0,x_0])$ satisfying $0 \leq g(x) \leq b(x)$ for $0 \leq x \leq x_0$, in which $s(x)$ is a smooth function that is positive on $[0,x_0)$ and decreases to $s(x_0) = 0$.

\subsubsection{Existence of allowed maximizers}
\label{sec:Proof_of_optimal_threshold_controls__Existence_of_allowed_maximizers}

\begin{proposition}
There exist maximizers $\criticalpoint{f} \in \mathcal{F}$ that maximize $R(f)$.
\end{proposition}

\myProof{
We will use several notions from the theory of Hilbert spaces and Lebesgue integration on the line. We consider maximization of $J(g)$ in \refEquation{eqn:Definition_of_Jg} over all measurable $g$ with $0 \leq g(x) \leq b(x)$ for a.e.\ $x \in [0,x_0]$. 

For any $g \in L^1( [0,x_0] )$, the Lebesgue points of $g$, i.e., all $x_1 \in (0,x_0)$ such that
\begin{equation}
\lim_{\varepsilon \downarrow 0} \frac{1}{2\varepsilon} \int_{-\varepsilon}^{\varepsilon} g(x_1+x) \d{x}
\label{eqn:Lebesgue_point_condition}
\end{equation} 
exists, is a subset of $[0,x_0]$ whose complement is a null set, and the limit in \refEquation{eqn:Lebesgue_point_condition} agrees with $g(x_1)$ for a.e.\ $x_1 \in [0,x_0]$, \cite{teschl_topics_2014}.

The set of allowed functions $g$ is a closed and bounded set of the separable Hilbert space $L^2([0,x_0])$, and the functional $J(g)$ is bounded on this set. Hence, we can find a sequence of candidates $\{ g_n \}_{ n \in \naturalNumbersZero }$ of allowed $g_n$, such that
\begin{equation}
\lim_{n \rightarrow \infty} J(g_n) = \sup_{ \textrm{allowed } g } \{ J(g) \} < \infty. 
\end{equation}
We can subsequently find a subsequence $\{ h_k \}_{ k \in \naturalNumbersZero } = \{ g_{n_k} \}_{ k \in \naturalNumbersZero }$ such that $h_k$ converges weakly to an $h \in L^2([0,x_0])$, \cite{rudin_real_1987}. Then
\begin{equation}
\sup_{ \textrm{allowed } g } \{ J(g) \}
= \lim_{k \rightarrow \infty} J(h_k) 
= \lim_{k \rightarrow \infty} \frac{ \int_0^\infty \bigl( \frac{r(x)}{A} - \frac{1}{B} \bigr) h_k(x) \d{x} }{ 1 + \frac{1}{B} \int_0^\infty h_k(x) \d{x} }
\overset{\mathrm{(i)}}= \frac{ \int_0^\infty \bigl( \frac{r(x)}{A} - \frac{1}{B} \bigr) h(x) \d{x} }{ 1 + \frac{1}{B} \int_0^\infty h(x) \d{x} }
= J(h),
\label{eqn:Candidate_h_is_maximizer}
\end{equation}
where (i) follows from weak convergence. We now only need to show that $h$ is allowed. We have for any $\varepsilon > 0$ and any $x_1 \in (0,x_0)$ by weak convergence that
\begin{equation}
\frac{1}{2\varepsilon} \int_{-\varepsilon}^{\varepsilon} h(x_1 + x) \d{x} 
= \lim_{k \rightarrow \infty} \frac{1}{2\varepsilon} \int_{-\varepsilon}^{\varepsilon} h_k(x_1 + x) \d{x} 
\in [0,b(x_1)],
\end{equation}
since all $h_k$ are allowed. Hence for all Lebesgue points $x_1 \in (0,x_0)$ of $h$ we have
\begin{equation}
\lim_{\varepsilon \downarrow 0} \frac{1}{2\varepsilon} \int_{-\varepsilon}^{\varepsilon} h(x_1 + x) \d{x} 
\in [0,b(x_1)],
\end{equation}
and so $0 \leq h(x_1) \leq b(x_1)$ for a.e. $x_1 \in [0,x_0]$. This, together with \refEquation{eqn:Candidate_h_is_maximizer} shows that $h$ is an allowed maximizer.
}

\subsubsection{Necessary condition for maximizers}
\label{sec:Proof_of_optimal_threshold_controls__Necessary_condition_for_maximizers}

\begin{proposition}
\label{prop:Necessary_condition_for_maximizers}
For any maximizer $\criticalpoint{f} \in \mathcal{F}$, $f(x) = 1$ if $r(x) > R(x)$, and $f(x) = 0$ if $r(x) < R(x)$. 
\end{proposition}

\myProof{
Let $g$ be an allowed maximizer of $J(g)$. We shall equivalently show that for any Lebesgue point $x_1 \in (0,x_0)$ of $g$,
\begin{equation}
\frac{ s(x_1) ( 1 + \int_0^{x_0} S g(x) \d{x} ) }{ S \int_0^{x_0} s(x) g(x) \d{x} } > 1 \Rightarrow g(x_1) = b(x_1),
%\quad \textrm{and}
\quad
\frac{ s(x_1) ( 1 + \int_0^{x_0} S g(x) \d{x} ) }{ S \int_0^{x_0} s(x) g(x) \d{x} } < 1 \Rightarrow g(x_1) = 0
\label{eqn:Equivalent_relations}
\end{equation}

Let $x_1 \in (0,x_0)$ be any Lebesgue point of $g$ and assume that
\begin{equation}
\frac{ s(x_1) ( 1 + \int_0^{x_0} S g(x) \d{x} ) }{ S \int_0^{x_0} s(x) g(x) \d{x} } > 1. \label{eqn:Assumption_on_ratio_of_s_integrals}
\end{equation}
Suppose that $g(x_1) < b(x_1)$. We shall derive a contradiction. Let $\varepsilon_0 > 0$ be small enough so that
\begin{equation}
\frac{1}{2} ( g(x_1) + b(x_1) ) \leq \min_{ x_1 - \varepsilon_0 \leq x \leq x_1 + \varepsilon_0 } \{ b(x) \}.
\label{eqn:Condition_of_varepsilon}
\end{equation}
Along with $g$, consider for $0 < \varepsilon \leq \varepsilon_0$ the function
\begin{equation}
g_{\varepsilon}(x) 
= 
\begin{cases}
g(x), & x \not\in [ x_1 - \varepsilon, x_1 + \varepsilon ], \\
\frac{1}{2} (  g(x_1) + b(x_1) ), & x \in [ x_1 - \varepsilon, x_1 + \varepsilon ]. \\
\end{cases}
\end{equation}
This $g_{\varepsilon}$ is allowed by \refEquation{eqn:Condition_of_varepsilon}. Write $J(g)$ as
\begin{equation}
J(g) 
= \frac{ C(\varepsilon) + I_s(\varepsilon; g) }{ D(\varepsilon) + I_S(\varepsilon; g) },
\end{equation}
where
\begin{equation}
C(\varepsilon) 
= \Biggl( \int_0^{y - \varepsilon} + \int_{y + \varepsilon}^{x_0} \Biggr) s(x) g(x) \d{x},
\quad
D(\varepsilon)
= 1 + \Biggl( \int_0^{y - \varepsilon} + \int_{y + \varepsilon}^{x_0} \Biggr) S g(x) \d{x},
\end{equation}
and
\begin{equation}
I_s(\varepsilon;g) 
= \int_{x_1-\varepsilon}^{x_1+\varepsilon} s(x) g(x) \d{x},
\quad
I_S(\varepsilon;g) 
= \int_{x_1-\varepsilon}^{x_1+\varepsilon} S g(x) \d{x}.
\label{eqn:Definition_of_Is_and_IS}
\end{equation}
We can do a similar thing with $J(g_{\varepsilon})$, using the same numbers $C(\varepsilon)$ and $D(\varepsilon)$ and $g$ replaced by $g_\varepsilon$ in \refEquation{eqn:Definition_of_Is_and_IS}. We compute
\begin{equation}
J(g_\varepsilon) - J(g)
= \frac{ ( C(\varepsilon) + I_s(\varepsilon; g_\varepsilon) ) ( D(\varepsilon) + I_S(\varepsilon;g) ) - ( C(\varepsilon) + I_s(\varepsilon; g) ) ( D(\varepsilon) + I_S(\varepsilon;g_\varepsilon) ) }{ ( D(\varepsilon) + I_S(\varepsilon;g) ) ( D(\varepsilon) + I_S(\varepsilon;g_\varepsilon) ) },
\label{eqn:Lebesgue_difference_of_Jg}
\end{equation}
in which the numerator $N(g_\varepsilon;g)$ of the fraction at the right-hand side of \refEquation{eqn:Lebesgue_difference_of_Jg} can be written as
\begin{align}
N(g_\varepsilon;g)
= & 
( I_s(\varepsilon; g_\varepsilon) - I_s(\varepsilon; g) ) D(\varepsilon) - ( I_S(\varepsilon; g_\varepsilon) - I_S(\varepsilon; g) ) C(\varepsilon) 
\nonumber \\ &
+ I_s(\varepsilon; g_\varepsilon) I_S(\varepsilon; g) - I_s(\varepsilon; g) I_S(\varepsilon; g_\varepsilon).
\end{align}
Since $x_1$ is a Lebesgue point of $g$, we have as $\varepsilon \downarrow 0$
\begin{gather}
\frac{1}{2\varepsilon} I_s(\varepsilon; g_\varepsilon) 
\rightarrow \frac{1}{2} s(x_1) ( g(x_1) + b(x_1) ),
\quad
\frac{1}{2\varepsilon} I_s(\varepsilon; g) 
\rightarrow s(x_1) g(x_1), \\
\frac{1}{2\varepsilon} I_S(\varepsilon; g_\varepsilon) 
\rightarrow \frac{1}{2} S ( g(x_1) + b(x_1) ),
\quad
\frac{1}{2\varepsilon} I_S(\varepsilon; g) 
\rightarrow S g(x_1),
\end{gather}
while also
\begin{equation}
C(\varepsilon) 
\rightarrow \int_0^{x_0} s(x) g(x) \d{x},
\quad
D(\varepsilon) 
\rightarrow 1 + \int_0^{x_0} S g(x) \d{x}.
\end{equation}
Therefore,
\begin{equation}
\lim_{\varepsilon \downarrow 0} N(g_{\varepsilon},g)
= \frac{1}{2} ( b(x_1) - g(x_1) ) \Bigl( s(x_1) \Bigl( 1 + \int_0^{x_0} S g(x) \d{x} \Bigr) - S \int_0^{x_0} s(x) g(x) \d{x} \Bigr) > 0
\end{equation}
by assumption \refEquation{eqn:Assumption_on_ratio_of_s_integrals}. Then $J(g_\varepsilon) - J(g) > 0$ when $\varepsilon$ is sufficiently small, contradicting maximality of $J(g)$. Hence, we have proven the first relation in \refEquation{eqn:Equivalent_relations}. The proof of the second relation is similar.
}

\subsubsection{Characterization of maximizers}
\label{sec:Proof_of_optimal_threshold_controls__Characterization_of_maximizers}

\refProposition{prop:Necessary_condition_for_maximizers} does not exclude the possibility that a maximizer alternates between $0$ and $1$. \refProposition{prop:Characterization_of_maximizers} solves this problem by excluding the pathological candidates.

\begin{proposition}
\label{prop:Characterization_of_maximizers}
The quantity 
\begin{equation}
R(f; \eta) 
= \frac{ A + \int_0^{\eta} r(x) f(x) \e{-\gamma x} \d{x} }{ B + \int_0^\eta f(x) \e{-\gamma x} \d{x} }.
\end{equation}
is uniquely maximized by 
\begin{equation}
f(x) = \indicator{ 0 \leq x \leq \criticalpoint{\eta} },
\end{equation}
with $\criticalpoint{\eta}$ a solution of the equation $r(\eta) = R_{\mathrm{T}}(\eta)$, apart from null functions and its value at any solution of $r(\eta) = R_{\mathrm{T}}(\eta)$.
\end{proposition}

\myProof{
Assume that $g$ is a maximizer, and consider the continuous, decreasing function
\begin{equation}
t_g(x_1) 
= s(x_1) \bigl( 1 + \int_0^{x_0} S g(x) \d{x} \bigr) - S \int_0^{x_0} s(x) g(x) \d{x},
\end{equation}
which is positive at $x_1 = 0$ and negative (because $g \neq 0$) at $x_1 = x_0$ since $s$ is decreasing with $s(0) > 0 = s(x_0)$. Let $x_{2,g}$, $x_{3,g}$ be such that $0 < x_{2,g} \leq x_{3,g} < x_0$ and 
\begin{equation}
t_g(x_1) 
= 
\begin{cases}
> 0, & 0 \leq x_1 < x_{2,g}, \\
= 0, & x_{2,g} \leq x_1 \leq x_{3,g}, \\
< 0, & x_{3,g} < x_1 \leq x_0. \\
\end{cases}
\end{equation}
Note that $x_{2,g} = x_{3,g}$ when $s$ is strictly decreasing on $[0,x_0]$, and that $s'(x) = 0$ for $x \in [x_{2,g}, x_{3,g} ]$ when $x_{2,g} < x_{3,g}$. According to \cite{janssen_scaled_2013}, we have
\begin{equation}
g(x_1) = b(x_1), 
\quad \textrm{a.e. } x_1 \in [0, x_{2,g}],
\quad \textrm{and} \quad
g(x_1) = 0, 
\quad \textrm{a.e. } x_1 \in [x_{3,g},x_0]. 
\label{eqn:Shape_of_g_between_0_x2g_x3g_and_x0}
\end{equation} 

For an allowed $h \neq 0$, consider the continuous function
\begin{equation}
J(h;x_1) 
= \frac{ \int_0^{x_1} s(x) h(x) \d{x} }{ 1 + \int_0^{x_1} S h(x) \d{x} }, 
\quad 0 \leq x_1 \leq x_0.
\label{eqn:Definition_of_Jhx1}
\end{equation}
We differentiate $J(h;x_1)$ with respect to $x_1$, where we use the fact that for any $k \in L^1([0,x_0])$,
\begin{equation}
\frac{\d{}}{\d{x_1}} \Bigl[ \int_0^{x_1} k(x) \d{x} \Bigr]
= k(x_1),
\quad \textrm{a.e. } x_1 \in [0,x_0].
\end{equation}
Thus we get for a.e.\ $x_1$ that
\begin{equation}
\frac{\d{}}{\d{x_1}} \Bigl[ J(h; x_1) \Bigr]
= \frac{N_h(x_1)}{D_h(x_1)},
\label{eqn:Derivative_of_Jhx1}
\end{equation}
where $D_h(x_1) = ( 1 + \int_0^{x_0} S h(x) \d{x} )^2$, and
\begin{equation}
N_h(x_1) = h(x_1) M_h(x_1)
\label{eqn:Numerator_of_derivative_of_Jhx1}
\end{equation}
with
\begin{equation}
M_h(x_1) 
= s(x_1) \bigl( 1 + \int_0^{x_1} S h(x) \d{x} \bigr) - S \int_0^{x_1} s(x) h(x) \d{x}.
\end{equation}
Now $M_h(x_1)$ is a continuous function of $x_1 \in [0,x_0]$ with $M_h(x_0) < 0 < M_h(0)$ since $s(x_0) = 0 < s(0)$ and $h \neq 0$. Furthermore, $M_h(x_1)$ is differentiable at a.e.\ $x_1$, and one computes for a.e.\ $x_1$,
\begin{equation}
\frac{\d{}}{\d{x_1}} \Bigl[ M_h(x_1) \Bigr]
= s'(x_1) \bigl( 1 + \int_0^{x_1} S h(x) \d{x} \bigr).
\label{eqn:Derivative_of_Mh}
\end{equation}
Since $s$ is decreasing, the right-hand side of \refEquation{eqn:Derivative_of_Mh} is nonpositive for all $x_1$ and negative for all $x_1$ with $s'(x_1) < 0$. 

Now let $g$ be a maximizer, and consider first the case that $s(x)$ is strictly decreasing. Then $x_{2,g} = x_{3,g}$ in \refEquation{eqn:Shape_of_g_between_0_x2g_x3g_and_x0}. Next consider $h = b$ in \refEquation{eqn:Definition_of_Jhx1} and further. It follows from \refEquation{eqn:Derivative_of_Mh} that $M_{b}$ is strictly decreasing on $[0,x_0]$, and so $M_{b}$ has a unique zero $\hat{x}$ on $[0,x_0]$. Therefore, by \refEquation{eqn:Derivative_of_Jhx1} and \refEquation{eqn:Numerator_of_derivative_of_Jhx1}, $J(b;x_1)$ has a unique maximum at $x_1 = \hat{x}$. Then, from \refEquation{eqn:Shape_of_g_between_0_x2g_x3g_and_x0} and maximality of $g$, $x_{2,g} = \hat{x} = x_{3,g}$. Hence, $J$ is uniquely maximized by 
\begin{equation}
g(x_1) 
= b(x_1) \indicator{ x_1 \in [0,\hat{x}] },
\label{eqn:Shape_of_gx1}
\end{equation}
apart from null functions, with $\hat{x}$ the unique solution $y$ of the equation
\begin{equation}
s(y) \bigl( 1 + \int_0^y S b(x) \d{x} \bigr) - S \int_0^y s(x) b(x) \d{x} = 0.
\label{eqn:Equation_for_y}
\end{equation}
This handles the case that $s$ is strictly decreasing.

When $s'$ may vanish, we have to argue more carefully. In the case that $x_{2,g} = x_{3,g}$, we can proceed as earlier, with \refEquation{eqn:Shape_of_gx1} emerging as maximizer and $x_{2,g} = y = x_{3,g}$. So assume we have a maximizer $g$ with $x_{2,g} < x_{3,g}$, and consider $h = g$ in \refEquation{eqn:Definition_of_Jhx1} and further. We have that $J(h=g; x_1)$ is constant in $x_1 \in [ x_{3,g}, x_0 ]$. Furthermore, from $s'(x_1) = 0$ for $x_1 \in [x_{2,g}, x_{3,g}]$ and \refEquation{eqn:Numerator_of_derivative_of_Jhx1}, we see that $J(h=g;x_1)$ is constant in $x_1 \in [x_{2,g}, x_{3,g}]$ as well. This constant value equals $J(g)$, and is equal to $J( b \indicator{ x_1 \in [0,x_{2,g}] } )$ since, due to \refEquation{eqn:Numerator_of_derivative_of_Jhx1}, we have $J(g;\cdot)=J(\bar{g};\cdot)$ when $g = \bar{g}$ a.e.\ outside $[x_{2,g}, x_{3,g}]$. We are then again in the previous situation, and the solutions $y$ of \refEquation{eqn:Equation_for_y} form now a whole interval $[y_2, y_3]$. The maximizers are again unique, apart from their values for $x_1 \in [y_2,y_3]$ that can be chosen arbitrarily between $0$ and $b(x_1)$.
}

\section{Conclusions and future perspectives}

The \gls{QED} regime has gained tremendous popularity in the operations management literature, because it describes how large-scale service operations can achieve high system utilization while simultaneously maintaining short delays. Operating a system in the \gls{QED} regime typically entails hiring a number of servers according to the square-root staffing rule  $s=\lambda/\mu+\gamma \sqrt{\lambda/\mu}$, and has the added benefit that limiting performance measures can be described by elementary functions of just the one parameter $\gamma$. Through the square-root staffing rule, $\gamma$ determines a hedge against variability or overcapacity, which is of the order of the natural fluctuations of the demand per time unit when the system operates in the \gls{QED} regime. Classical problems of dimensioning large-scale systems in the \gls{QED} regime can then be solved by optimizing objective functions solely dependent on $\gamma$.

Our paper adds a revenue maximization framework that complies with the classical dimensioning of \gls{QED} systems by constructing scalable admission controls and revenue structures that remain meaningful in the \gls{QED} regime (\refProposition{prop:QED_limit_of_revenue__Continuous_admission_controls}). 
As we have proven, our revenue framework naturally leads to an optimal control that
bars customers from entering when the queue length of delayed customers exceeds the threshold $\eta^{\rm opt}\sqrt{s}$, provided that $\eta^{\rm opt}$ satisfies a fundamental \emph{threshold equation} (\refProposition{prop:Threshold_control_is_optimal_form_of_control}). A detailed study of this threshold equation made it possible to characterize $\eta^{\rm opt}$ in terms of exact expressions, bounds, and asymptotic expansions. The weak assumptions made throughout this paper allow for application to a rich class of revenue structures, and an interesting direction for future work would therefore be the construction of realistic revenue structures based on specific case studies, expert opinions, or calibration to financial data.  

Let us finally discuss the fascinating interplay between the parameters $\gamma$ and $\eta$, which suggest that they act as communicating yet incomparable vessels. The optimal threshold $\eta^{\rm opt}$ increases with the overcapacity $\gamma$. Since more overcapacity roughly means fewer customers per server, and a larger threshold means more customers per server, we see that the optimization of revenues over the pair $(\gamma,\eta)$ gives rise to an intricate two-dimensional framework in which the two parameters have radically different yet persistent effects in the \gls{QED} regime. At the process level, the $\gamma$ acts as a negative drift in the entire state space, while the $\eta$ only interferes at the upper limit of the state space. Hence, while in this paper we have treated $\gamma$ as given, and mostly focused on the behavior of the new parameter $\eta$, our framework paves the way for two-dimensional joint staffing and admission control problems. Gaining a deeper understanding of this interplay, and in relation to specific revenue structures, is a promising direction for future research.

\myAppendices{

\section{Limiting behavior of long-term \gls{QED} revenue}
\label{sec:Limiting_behavior_of_long_term_QED_revenue}

With $r_s(k) = r( (k-s)/\sqrt{s} )$ as in \refEquation{formr} and $\pi_s(k) = \lim_{t \to \infty} \probability{ Q_s(t) = k }$, \refEquation{eqn:Equilibrium_distribution}, where $p_s$ and $f$ are related as in \refEquation{formf}, we compute for $\rho = 1 - \gamma / \sqrt{s} > 0$,
\begin{equation}
\sum_{k=0}^\infty r_s(k) \pi_s(k)
= \frac{ \sum_{k=0}^s r\bigl( \frac{k-s}{\sqrt{s}} \bigr) \frac{(s\rho)^k}{k!} + \frac{(s\rho)^s}{s!} \sum_{k=s+1}^\infty r\bigl( \frac{k-s}{\sqrt{s}} \bigr) \rho^{k-s} f\bigl( \frac{k-s}{\sqrt{s}} \bigr) }{ \sum_{k=0}^s \frac{(s\rho)^k}{k!} + \frac{(s\rho)^s}{s!} \sum_{k=s+1}^\infty \rho^{k-s} f\bigl( \frac{k-s}{\sqrt{s}} \bigr) }.
\end{equation}
Dividing by the factor $(s\rho)^s / s!$, we obtain
\begin{equation}
\sum_{k=0}^\infty r_s(k) \pi_s(k)
= \frac{W_s^{\mathrm{L}}(\rho) + W_s^{\mathrm{R}}(\rho)}{B_s^{-1}(\rho) + F_s(\rho)}.
\end{equation}
Here,
\begin{equation}
B_s(\rho) = \frac{\frac{(s\rho)^s}{s!}}{\sum_{k=0}^s \frac{(s\rho)^k}{k!}}
\end{equation}
is  the Erlang B formula,
\begin{equation}
F_s(\rho)
= \sum_{n=0}^\infty \rho^{n+1} f\Bigl( \frac{n+1}{\sqrt{s}} \Bigr)
\end{equation}
as in \refEquation{eqn:Fs}, and
\begin{gather}
W_s^{\mathrm{L}}(\rho) 
= \sum_{k=0}^s r\Bigl( \frac{k-s}{\sqrt{s}} \Bigr) \frac{s!(s\rho)^{k-s}}{k!},
\label{eqn:Definition__WsL}
\\
W_s^{\mathrm{R}}(\rho)
= \sum_{n=0}^\infty r\Bigl( \frac{n+1}{\sqrt{s}} \Bigr) \rho^{n+1} f\Bigl( \frac{n+1}{\sqrt{s}} \Bigr).
\label{eqn:Definition__WsR}
\end{gather}
with superscripts $\mathrm{L}$ and $\mathrm{R}$ referring to the left-hand part $k=0, 1, \ldots, s$ and right-hand part $k=s+1, s+2, \ldots$ of the summation range, respectively.

From Jagerman's asymptotic results for Erlang B, there is the approximation \cite[Theorem~14]{jagerman_properties_1974}
\begin{equation}
B_s^{-1}(\rho) 
= \sqrt{s} \psi(\gamma) + \chi(\gamma) + \bigObig{ \frac{1}{\sqrt{s}} } 
\label{eqn:Asymptotic_expansion_of_Bs_inverse}
\end{equation}
with $\psi(\gamma) = \Phi(\gamma) / \phi(\gamma)$ and $\chi(\gamma)$ expressible in terms of $\phi$ and $\Phi$ as well. For $F_s(\rho)$ there is the approximation \cite[Theorem~4.2]{janssen_scaled_2013},
\begin{equation}
F_s(\rho) 
= \sqrt{s} \mathcal{L}(\gamma) + \mathcal{M}(\gamma) + \bigObig{\frac{1}{\sqrt{s}}},
\label{eqn:Asymptotic_expansion_of_Fs}
\end{equation}
with $\mathcal{L}(\gamma) = \int_0^\infty f(x) \exp{(- \gamma x)} \d{x}$ and $\mathcal{M}(\gamma)$ expressible in terms of $\mathcal{L}'(\gamma)$. We aim at similar approximations for $W_s^{\mathrm{L}}(\rho)$ and $W_s^{\mathrm{R}}(\rho)$ in \refEquation{eqn:Definition__WsL}, \refEquation{eqn:Definition__WsR}.

We start by considering $W_s^{\mathrm{R}}(\rho)$ for the case that $r$ and its first two derivatives are continuous and bounded in the two following situations:
\begin{itemize}
\item[$(\mathrm{i.})$] $f$ is smooth; $f(y) \exp{(-\gamma y)}$ and its first two derivatives are exponentially small as $y \rightarrow \infty$.
\item[$(\mathrm{ii.})$] $f = \indicator{ x \in [0,\eta] }$ with $\eta > 0$.
\end{itemize}

\subsection{Asymptotics of $W_s^{\mathrm{R}}(\rho)$}

In the series expression for $W_s^{\mathrm{R}}(\rho)$, we have 
\begin{equation}
\rho^{n+1} 
= \Bigl( 1 - \frac{\gamma}{\sqrt{s}} \Bigr)^{n+1}
= \e{-\frac{(n+1)\gamma_s}{\sqrt{s}}}
\end{equation}
with
\begin{equation}
\gamma_s 
= - \sqrt{s} \ln{ \Bigl( 1 - \frac{\gamma}{\sqrt{s}} \Bigr) } 
= \gamma + \frac{\gamma^2}{2\sqrt{s}} + \ldots > \gamma.
\label{eqn:Expansion_of_gammas}
\end{equation}
Hence, the conditions in case $(\mathrm{i.})$ are also valid when using $\gamma_s$ instead of $\gamma$.

We obtain the following result.
\begin{lemma}
\label{lem:Expansions_of_WsR_for_case_I_and_II}
For case $(\mathrm{i.})$ it holds that
\begin{equation}
W_s^{\mathrm{R}}(\rho)
= \sqrt{s} \int_0^\infty \e{-\gamma_s y} r(y) f(y) \d{y} - \frac{1}{2} r(0) f(0) + \bigObig{\frac{1}{\sqrt{s}}}.
\label{eqn:WsR_expansion__Case_I}
\end{equation}
For case $(\mathrm{ii.})$ it holds that
\begin{equation}
W_s^{\mathrm{R}}(\rho)
= \sqrt{s} \int_0^\eta \e{-\gamma_s y} r(y) \d{y} - \frac{1}{2} r(0) + \Bigl( \lfloor \eta \sqrt{s} \rfloor - \bigl( \eta\sqrt{s} - \frac{1}{2} \bigr) \Bigr) \e{-\gamma_s \eta} r(\eta) + \bigObig{\frac{1}{\sqrt{s}}}.
\label{eqn:WsR_expansion__Case_II}
\end{equation}
\end{lemma}

\myProof{
We use \gls{EM}-summation as in \cite[Appendix C]{janssen_scaled_2013}, first instance in \cite[(C.1)]{janssen_scaled_2013}, case $m=1$, with the function
\begin{equation}
h(x) 
= g\Bigl( \frac{x+\frac{1}{2}}{\sqrt{s}} \Bigr),
\quad x \geq 0,
\quad \textrm{and} \quad
g(y) 
= \e{-\gamma_s y} r(y)f(y), 
\quad y \geq 0,
\end{equation}
using a finite summation range $n=0,1, \ldots, N$, where we take $N=s$ in case $(\mathrm{i.})$ and $N = \lfloor \eta \sqrt{s} - 3/2 \rfloor$ in case $(\mathrm{ii.})$ In both cases, we have by smoothness of $h$ on the range $[0,N+1]$ that
\begin{equation}
\sum_{n=0}^N h\bigl( n + \tfrac{1}{2} \bigr) 
= \int_0^{N+1} h(x) \d{x} + \tfrac{1}{2} B_2\bigl( \tfrac{1}{2} \bigr) \bigl(  h^{(1)}(N+1) - h^{(1)}(0) \bigr) + R,
\label{eqn:Summation_over_h}
\end{equation}
where $|R| \leq \frac{1}{2} B_2 \int_0^{N+1} | h^{(2)}(x) | \d{x}$. Due to our assumptions, it holds in both cases that
\begin{equation}
\tfrac{1}{2} B_2\bigl( \tfrac{1}{2} \bigr) \bigl(  h^{(1)}(N+1) - h^{(1)}(0) \bigr) + R
= \bigObig{\frac{1}{\sqrt{s}}}.
\end{equation}
In case $(\mathrm{i.})$, the left-hand side of \refEquation{eqn:Summation_over_h} equals $W_s^{\mathrm{R}}(\rho)$, apart from an error that is exponentially small as $s \rightarrow \infty$. In case $(\mathrm{ii.})$, the left-hand side of \refEquation{eqn:Summation_over_h} and $W_s^{\mathrm{R}}(\rho)$ are related according to
\begin{equation}
W_s^{\mathrm{R}}(\rho)
= \sum_{n=0}^N h\bigl( n + \tfrac{1}{2} \bigr) + g\Bigl( \frac{ \lfloor \eta \sqrt{s} \rfloor }{\sqrt{s}} \Bigr) \bigl( \lfloor \eta \sqrt{s} \rfloor - \bigl\lfloor \eta\sqrt{s} - \tfrac{1}{2} \bigr\rfloor \bigr).
\label{eqn:Relation_between_LHS_and_WsR}
\end{equation}
The second term at the right-hand side of \refEquation{eqn:Relation_between_LHS_and_WsR} equals $0$ or $g( \lfloor \eta \sqrt{s} \rfloor / \sqrt{s} )$ accordingly as $\eta \sqrt{s} - \lfloor \eta \sqrt{s} \rfloor \geq$ or $< \frac{1}{2}$, i.e., accordingly as $N+1 = \lfloor \eta \sqrt{s} \rfloor$ or $\lfloor \eta \sqrt{s} \rfloor - 1$. Next, by smoothness of $h$ and $g$ on the relevant ranges, we have
\begin{equation}
\int_0^{N+1} h(x) \d{x} 
= \sqrt{s} \int_{\frac{1}{2\sqrt{s}}}^{\frac{N+3/2}{\sqrt{s}}} g(y) \d{y}
= \sqrt{s} \int_0^{\frac{N+3/2}{\sqrt{s}}} g(y) \d{y} - \frac{1}{2} g(0) + \bigObig{\frac{1}{\sqrt{s}}}.
\end{equation}
In case $(\mathrm{i.})$, we have that $\int_{(N+3/2)/\sqrt{s}}^\infty g(y) \d{y}$ is exponentially small as $s \rightarrow \infty$, since $N=s$, and this yields \refEquation{eqn:WsR_expansion__Case_I}. In case $(\mathrm{ii.})$, we have
\begin{align}
\int_0^{\frac{N+3/2}{\sqrt{s}}} g(y) \d{y} - \int_0^\eta g(y) \d{y} 
&
= \int_\eta^{\frac{N+3/2}{\sqrt{s}}} g(y) \d{y}
\nonumber \\ &
= \Bigl( \frac{N+3/2}{\sqrt{s}} - \eta \Bigr) g\Bigl( \frac{ \lfloor \eta \sqrt{s} \rfloor}{\sqrt{s}} \Bigr) + \bigObig{\frac{1}{s}}
\nonumber \\ &
= \frac{1}{\sqrt{s}} \bigl( \bigl\lfloor \eta \sqrt{s} - \tfrac{1}{2} \bigr\rfloor - \bigl( \eta \sqrt{s} - \tfrac{1}{2} \bigr) \bigr) g\Bigl( \frac{\lfloor \eta \sqrt{s} \rfloor}{\sqrt{s}} \Bigr) + \bigObig{\frac{1}{s}},
\end{align}
and with \refEquation{eqn:Relation_between_LHS_and_WsR}, this yields \refEquation{eqn:WsR_expansion__Case_II}. This completes the proof.
}

We denote for both case $(\mathrm{i.})$ and $(\mathrm{ii.})$
\begin{equation}
\mathcal{L}_{rf}(\delta) 
= \int_0^\infty \e{-\delta y} r(y) f(y) \d{y}
\label{eqn:Definition__Lrf}
\end{equation}
with $\delta \in \realNumbers$ such that the integral of the right-hand side of \refEquation{eqn:Definition__Lrf} converges absolutely. From \refEquation{eqn:Expansion_of_gammas} it is seen that, with the prime $'$ denoting differentiation,
\begin{equation}
\mathcal{L}_{rf}(\gamma_s)
= \mathcal{L}_{rf}(\gamma) + \frac{\gamma^2}{2\sqrt{s}} \mathcal{L}_{rf}'(\gamma) + \bigObig{\frac{1}{s}}.
\end{equation}
Thus we get from \refLemma{lem:Expansions_of_WsR_for_case_I_and_II} the following result.

\begin{proposition}
\label{prop:Asymptotic_expansion_of_WsR}
For case $(\mathrm{i.})$ it holds that
\begin{equation}
W_s^{\mathrm{R}}(\rho)
= \sqrt{s} \mathcal{L}_{rf}(\gamma) + \tfrac{1}{2} \gamma^2 \mathcal{L}_{rf}'(\gamma) - \tfrac{1}{2} r(0) f(0) + \bigObig{\frac{1}{\sqrt{s}}}.
\end{equation}
For case $(\mathrm{ii.})$ it holds that
\begin{equation}
W_s^{\mathrm{R}}(\rho)
= \sqrt{s} \mathcal{L}_{rf}(\gamma) + \tfrac{1}{2} \gamma^2 \mathcal{L}_{rf}'(\gamma) - \tfrac{1}{2} r(0) + \bigl( \lfloor \eta \sqrt{s} \rfloor - \bigl( \eta \sqrt{s} - \tfrac{1}{2} \bigr) \bigr) \e{-\gamma \eta} r(\eta) + \bigObig{\frac{1}{\sqrt{s}}}.
\end{equation}
\end{proposition}

\subsection{Asymptotics of $W_s^{\mathrm{L}}(\rho)$}

We next consider $W_s^{\mathrm{L}}(\rho)$ for the case that $r: (-\infty,0] \rightarrow \realNumbers$ has bounded and continuous derivatives up to order $2$. Using a change of variables, we write
\begin{equation}
W_s^{\mathrm{L}}(\rho)
= r(0) + \sum_{k=1}^s r\Bigl( \frac{-k}{\sqrt{s}} \Bigr) \frac{s!s^{-k}}{(s-k)!} \rho^{-k},
\label{eqn:Expression_for_WsL}
\end{equation}
and we again intend to apply \gls{EM}-summation to the series at the right-hand side of \refEquation{eqn:Expression_for_WsL}. We first present a bound and an approximation.

\begin{lemma}
We have for $|\gamma| / \sqrt{s} \leq \frac{1}{2}$ and $\rho = 1 - \gamma / \sqrt{s}$,
\begin{equation}
\frac{s!s^{-k}}{(s-k)!} \rho^{-k}
\leq \exp{ \Bigl( - \frac{k(k-1)}{2s} + \frac{\gamma k}{\sqrt{s}} + \frac{\gamma^2 k}{s} \Bigr) },
\quad k = 1, 2, \ldots, s,
\label{eqn:Exponential_bound_in_k_and_gamma_for_local_Poisson_ratio}
\end{equation}
and
\begin{equation}
\frac{s!s^{-k}}{(s-k)!} \rho^{-k}
= G_s\Bigl( \frac{k}{\sqrt{s}} \Bigr) \bigl( 1 + \bigObig{ \frac{1}{s} P_6\Bigl( \frac{k}{\sqrt{s}} \Bigr) } \bigr),
\quad k \leq s^{2/3},
\label{eqn:Exponential_bound_in_Gs_for_local_Poisson_ratio}
\end{equation}
where
\begin{equation}
G_s(y) 
= \e{ - \frac{1}{2} y^2 + \gamma y } \bigl( 1 - \frac{1}{6\sqrt{s}} y^3 + \frac{1}{2\sqrt{s}} (1 + \gamma^2) y \bigr),
\label{eqn:Definition_of_Gs}
\end{equation}
and $P_6(y)$ is a polynomial in $y$ of degree $6$ with coefficients bounded by $1$ (the constant implied by $\bigO{\cdot}$ depends on $\gamma$).
\end{lemma}

\myProof{
We have for $k=1,2,\ldots,s$ and $|\gamma| / \sqrt{s} \leq 1/2$, $\rho = 1 - \gamma / \sqrt{s}$,
\begin{align}
\frac{s!s^{-k}}{(s-k)!} \rho^{-k}
&
= \rho^{-k} \prod_{j=0}^{k-1} \Bigl( 1 - \frac{j}{s} \Bigr)
%\nonumber \\ &
= \exp{ \Bigl( \sum_{j=0}^{k-1} \ln{ \Bigl( 1 - \frac{j}{s} \Bigr) } - k \ln{ \Bigl( 1 - \frac{\gamma}{\sqrt{s}} \Bigr) } \Bigr) }
\nonumber \\ &
\leq \exp{ \Bigl( - \sum_{j=0}^{k-1} \frac{j}{s} + \frac{\gamma k}{\sqrt{s}} + \frac{\gamma^2 k}{s} \Bigr) }
%\nonumber \\ &
=  \exp{ \Bigl( - \frac{k(k-1)}{2s} + \frac{\gamma k}{\sqrt{s}} + \frac{\gamma^2 k}{s} \Bigr) },
\end{align}
where it has been used that $-\ln{(1-x)} \leq x + x^2$, $|x| \leq 1/2$. 

On the range $k \leq s^{2/3}$, we further expand
\begin{align}
\frac{s!s^{-k}}{(s-k)!} \rho^{-k}
&
= \exp{ \Bigl( - \sum_{j=0}^{k-1} \Bigl( \frac{j}{s} + \frac{j^2}{2s^2} + \bigObig{ \frac{j^3}{s^3} } \Bigr) + \frac{\gamma k}{\sqrt{s}} + \frac{\gamma^2 k}{s} + \bigObig{ \frac{k}{s^{3/2}} } \Bigr) }
\nonumber \\ &
= \exp{ \Bigl( - \frac{k(k-1)}{2s} - \frac{k(k-1)(2k-1)}{12s^2} + \bigObig{ \frac{k^4}{s^3} } + \frac{\gamma k}{\sqrt{s}} + \frac{\gamma^2 k}{s} + \bigObig{ \frac{k}{s^{3/2}} } \Bigr) }
\nonumber \\ &
= \exp{ \Bigl( - \frac{k^2}{2s} + \frac{\gamma k}{\sqrt{s}} - \frac{k^3}{6s^2} + \frac{1}{2} ( 1 + \gamma^2 ) \frac{k}{s} + \bigObig{ \frac{k}{s^{3/2}} + \frac{k^2}{s^2} + \frac{k^4}{s^3} } \Bigr) }
\end{align}
On the range $0 \leq k \leq s^{2/3}$ we have
\begin{equation}
\frac{k^3}{s^2}, \frac{k}{s}, \frac{k}{s^{3/2}}, \frac{k^2}{s^2}, \frac{k^4}{s^3} 
= \bigO{1}. 
\end{equation}
Hence, on the range $0 \leq k \leq s^{2/3}$,
\begin{align}
\frac{s!s^{-k}}{(s-k)!} \rho^{-k}
&
= \exp{ \Bigl( - \frac{k^2}{2s} + \frac{\gamma k}{\sqrt{s}} \Bigr) } \Bigl( 1 - \frac{k^3}{6s^2} + \frac{1}{2} ( 1 + \gamma^2 ) \frac{k}{s} + \bigObig{ \frac{k^6}{s^4} + \frac{k^4}{s^3} + \frac{k^2}{s^2} + \frac{k}{s^{3/2}} } \Bigr)
\nonumber \\ &
= G\Bigl( \frac{k}{\sqrt{s}} \Bigr) \Bigl( 1 + \bigObig{ \frac{1}{s} P_6\Bigl( \frac{k}{\sqrt{s}} \Bigr) } \Bigr),
\end{align}
where $P_6(y) = y^6 + y^4 + y^2 + y$.
}

\begin{proposition}
\label{prop:Asymptotic_expansion_of_WsL}
It holds that
\begin{equation}
W_s^{\mathrm{L}}(\rho)
= \sqrt{s} \int_{-\infty}^0 \e{-\frac{1}{2} y^2 - \gamma y} r(y) \d{y} + \frac{1}{2} r(0) + \int_{-\infty}^0 \e{-\frac{1}{2} y^2 - \gamma y } ( \frac{1}{6} y^3 - \frac{1}{2} ( 1 + \gamma^2 ) y ) r(y) \d{y} + \bigObig{ \frac{1}{\sqrt{s}} }. 
\end{equation}
\end{proposition}

\myProof{
With $v(y) = r(-y)$, we write
\begin{equation}
W_s^{\mathrm{L}}(\rho)
= r(0) + \sum_{n=0}^{s-1} v\Bigl( \frac{n+1}{\sqrt{s}} \Bigr) \frac{s!s^{-n-1}}{(s-n-1)!} \rho^{-n-1}.
\label{eqn:Expansion_of_WsL}
\end{equation}
By the assumptions on $r$ and the bound in \refEquation{eqn:Exponential_bound_in_k_and_gamma_for_local_Poisson_ratio}, the contribution of the terms in the series in \refEquation{eqn:Expansion_of_WsL} is $\bigO{ \exp{( - C s^{1/3} )} }$, $s \rightarrow \infty$, for any $C$ with $0 < C < 1/2$. On the range $n = 0, 1, \ldots, \lfloor s^{2/3} \rfloor - 1 =: N$, we can apply \refEquation{eqn:Exponential_bound_in_Gs_for_local_Poisson_ratio}, and so, with exponentially small error,
\begin{equation}
W_s^{\mathrm{L}}(\rho)
= r(0) + \sum_{n=0}^N v\Bigl( \frac{n+1}{\sqrt{s}} \Bigr) G_s\Bigl( \frac{n+1}{\sqrt{s}} \Bigr) \Bigl( 1 + \bigObig{ \frac{1}{s} P_6\Bigl( \frac{n+1}{\sqrt{s}} \Bigr) } \Bigr).
\end{equation}
By \gls{EM}-summation, as used in the proof of \refLemma{lem:Expansions_of_WsR_for_case_I_and_II} for the case $(\mathrm{i.})$ as considered there, we have
\begin{equation}
\sum_{n=0}^N v\Bigl( \frac{n+1}{\sqrt{s}} \Bigr) G_s\Bigl( \frac{n+1}{\sqrt{s}} \Bigr)
= \sqrt{s} \int_0^\infty v(y) G_s(y) \d{y} - \frac{1}{2} v(0) G_s(0) + \bigObig{\frac{1}{\sqrt{s}}},
\end{equation}
where we have extended the integration range $[0,(N+3/2)/\sqrt{s}]$ to $[0,\infty)$ at the expense of exponentially small error. Then the result follows on a change of the integration variable, noting that $v(y) = r(-y)$ and the definition of $G_s$ in \refEquation{eqn:Definition_of_Gs}, implying $G_s(0) = 1$.
}

The result of \refProposition{prop:QED_limit_of_revenue__Continuous_admission_controls} in the main text follows now from \refEquation{eqn:Asymptotic_expansion_of_Bs_inverse}, \refEquation{eqn:Asymptotic_expansion_of_Fs}, \refProposition{prop:Asymptotic_expansion_of_WsR} and \refProposition{prop:Asymptotic_expansion_of_WsL}, by considering leading terms only.

\section{Explicit solutions for exponential revenue}
\label{sec:Appendix__Explicit_solutions_for_exponential_revenue}

\begin{proposition}
\label{prop:Explicit_characterization_in_exponential_case_when_R0_is_approximately_one}
When $\varepsilon = 1-R_{\mathrm{T}}(0) > 0$ is sufficiently small, 
\begin{equation}
\criticalpoint{\eta} 
= - \frac{1}{\delta} \ln{ ( 1 - \sum_{l=1}^\infty a_l \varepsilon^l ) },
\label{eqn:Exponential_revenue__Series_solution}
\end{equation}
where
\begin{gather}
a_1 = 1,
\quad a_2 = \frac{1}{2} ( \alpha + \beta - 1 ), \label{eqn:Recursion_for_a_l__Initial_values}
\\
\quad a_{l+1} = \frac{1}{l+1} \Bigl( ( l \alpha + (l+1) \beta - 1 ) a_l + \beta \sum_{i=2}^{l-1} i a_i a_{l+1-i} \Bigr),
\quad l = 2, 3, \ldots,
\label{eqn:Recursion_for_a_l__Recursion}
\end{gather}
with $\beta = (1-\alpha)(1+1/(\gamma B))$ and the convention that $\sum_{i=2}^{l-1} = 0$ for $l=2$.
\end{proposition}

\myProof{
With $\varepsilon = 1 - R_{\mathrm{T}}(0)$ and $w = 1 - z$, we can write \refEquation{eqn:Threshold_equation_for_exponential_revenue_in_terms_of_z} as
\begin{equation}
H(w) 
= w + \frac{1}{\gamma B} ( w - \frac{1}{\alpha} ( 1 - (1-w)^\alpha ) ) 
= \varepsilon.
\label{eqn:Exponential_revenue__Definition_of_Hw}
\end{equation}
Note that 
\begin{equation}
H(w) 
= w + \frac{1}{\gamma B} ( \frac{1}{2} ( \alpha - 1 ) w^2 - \frac{1}{6} ( \alpha - 1) ( \alpha - 2 ) w^3 + \ldots ), 
\quad |w| < 1,
\end{equation}
and so there is indeed a (unique) solution
\begin{equation}
w(\varepsilon) 
= \varepsilon + \sum_{l=2}^\infty a_l \varepsilon^l
\label{eqn:Exponential_revenue__Power_series_for_w}
\end{equation}
of \refEquation{eqn:Exponential_revenue__Definition_of_Hw} when $|\varepsilon|$ is sufficiently small. To find the $a_l$ we let $D = 1 / ( \gamma B )$, and we write \refEquation{eqn:Exponential_revenue__Definition_of_Hw} as
\begin{equation}
(1+D) w - \frac{1}{\alpha} D + \frac{1}{\alpha} ( 1 - w )^\alpha 
= \varepsilon,
\quad 
w = w(\varepsilon).
\label{eqn:Exponential_revenue__Alternative_form_of_Hw}
\end{equation}
Differentiating \refEquation{eqn:Exponential_revenue__Alternative_form_of_Hw} with respect to $\varepsilon$, multiplying by $1 - w(\varepsilon)$, and eliminating $(1-w(\varepsilon))^\alpha$ using \refEquation{eqn:Exponential_revenue__Alternative_form_of_Hw} yields the equation
\begin{equation}
( 1 - \alpha \varepsilon - \beta w(\varepsilon) ) w'(\varepsilon) = 1 - w(\varepsilon).
\label{eqn:Exponential_revenue__Derivative_of_and_then_simplified_alternative_form_of_Hw}
\end{equation}
Inserting the power series \refEquation{eqn:Exponential_revenue__Power_series_for_w} for $w(\varepsilon)$ and $1 + \sum_{l=1}^\infty (l+1) a_{l+1} \varepsilon^l$ for $w'(\varepsilon)$ into \refEquation{eqn:Exponential_revenue__Derivative_of_and_then_simplified_alternative_form_of_Hw} gives
\begin{align}
&
1 - \alpha \varepsilon + \sum_{l=1}^\infty (l+1) a_{l+1} \varepsilon^l - \alpha \sum_{l=2}^\infty l a_l \varepsilon^l 
\nonumber \\ &
- \beta \varepsilon - \beta \sum_{l=2}^\infty l a_l \varepsilon^l - \beta \sum_{l=2}^\infty a_l \varepsilon^l - \beta \sum_{l=2}^\infty a_l \varepsilon^l \sum_{l=1}^\infty (l+1) a_{l+1} \varepsilon^l
\nonumber \\ &
= 1 - \varepsilon - \sum_{l=2}^\infty a_l \varepsilon^l.
\label{eqn:Exponential_revenue_DSAF_of_Hw_with_power_series}
\end{align}
Using that
\begin{equation}
\sum_{l=2}^\infty a_l \varepsilon^l \sum_{l=1}^\infty (l+1) a_{l+1} \varepsilon^l 
= \sum_{l=3}^\infty \Bigl( \sum_{i=2}^{l-1} i a_i a_{l+1-i} \Bigr) \varepsilon^l,
\end{equation}
it follows that $a_1, a_2, a_3, \ldots$ can be found recursively as in \refEquation{eqn:Recursion_for_a_l__Initial_values}--\refEquation{eqn:Recursion_for_a_l__Recursion}, by equating coefficients in \refEquation{eqn:Exponential_revenue_DSAF_of_Hw_with_power_series}. The result \refEquation{eqn:Exponential_revenue__Series_solution} then follows from $\criticalpoint{\eta} = - (1/\delta) \ln{z} = (1/\delta) \ln{(1-w)}$. The inequality $\beta < 0$ follows from the inequality $\gamma + \phi(\gamma) / \Phi(\gamma) > 0$, $\gamma \in \realNumbers$, given in \cite[Sec.~4]{avram_loss_2013}.
}

We consider next the cases $\alpha = -1$, $1/2$, and $2$ that allow for solving the threshold equation explicitly, and that illustrate \refProposition{prop:Explicit_characterization_in_exponential_case_when_R0_is_approximately_one}.

\begin{proposition}
\label{prop:Exponential_revenue__Square_root_solutions}
Let $t = - \gamma B / ( 1 + \gamma B ) >0$, and $\varepsilon = 1 - R_{\mathrm{T}}(0)$, for the cases $\mathrm{(i)}$ and $\mathrm{(ii)}$ below. The optimal threshold $\criticalpoint{\eta}$ is given as
\begin{equation}
\criticalpoint{\eta} 
= - \frac{1}{\delta} \ln{( 1 - w(\varepsilon) )},
\end{equation}
where $w(\varepsilon)$ is given by: \\
$\mathrm{(i)}$ $\alpha = -1$,
\begin{align}
w(\varepsilon)
&
= \frac{1}{2} t \Bigl( \sqrt{ (1+\varepsilon)^2 + \frac{4\varepsilon}{t}  } - 1 - \varepsilon \Bigr)
\nonumber \\ &
= \varepsilon - \frac{1}{2} t \sum_{k=2}^\infty (-1)^k \frac{ P_k \Bigl( 1 + \frac{2}{t} \Bigr) - P_{k-2} \Bigl( 1 + \frac{2}{t} \Bigr) }{2k-1} \varepsilon^k,
\quad |\varepsilon| < 1 + \frac{2}{t} - \sqrt{ \Bigl( 1 + \frac{2}{t} \Bigr)^2 - 1 },
\label{eqn:Exponential_revenue__Critical_w__Alpha_minus_1}
\end{align}
where $P_k$ is the Legendre polynomial of degree $k$,\\
$\mathrm{(ii)}$ $\alpha = 1/2$,
\begin{equation}
w(\varepsilon)
= \frac{2t}{1 + \gamma B} \Bigl( \sqrt{ 1 + \frac{\varepsilon}{t} } - 1 \Bigr) - t \varepsilon 
=  \varepsilon + \frac{2t}{1 + \gamma B} \sum_{k=2}^\infty \binom{1/2}{k} \Bigl( \frac{\varepsilon}{t} \Bigr)^k, 
\quad |\varepsilon| < t,
\label{eqn:Exponential_revenue__Critical_w__Alpha_half}
\end{equation}
$\mathrm{(iii)}$ $\alpha = 2$,
\begin{equation}
w(\varepsilon) 
= - \gamma B + \sqrt{ (\gamma B)^2 + 2 \gamma B \varepsilon }
= \varepsilon + \gamma B \sum_{k=2}^\infty \binom{1/2}{k} \Bigl( \frac{2\varepsilon}{\gamma B} \Bigr)^k,
\quad | \varepsilon | < \frac{1}{2} \gamma B.
\label{eqn:Exponential_revenue__Critical_w__Alpha_2}
\end{equation}
\end{proposition}

\myProof{
Case $\mathrm{(i)}$. When $\alpha = -1$, we can write the threshold equation as 
\begin{equation}
w^2 + t ( 1 + \varepsilon ) w = t \varepsilon.
\label{eqn:Exponential_revenue__Threshold_equation__Alpha_minus_1}
\end{equation}
From the two solutions 
\begin{equation}
w = - \frac{1}{2} t ( 1 + \varepsilon ) \pm \sqrt{ ( \frac{1}{2} t ( 1 + \varepsilon ) )^2 + t \varepsilon }
\end{equation}
of \refEquation{eqn:Exponential_revenue__Threshold_equation__Alpha_minus_1}, we take the one with the $+$ sign so as to get $w$ small and positive when $\varepsilon$ is small and positive. This gives $w(\varepsilon)$ as in the first line of \refEquation{eqn:Exponential_revenue__Critical_w__Alpha_minus_1}, the solution being analytic in the $\varepsilon$-range given in the second line of \refEquation{eqn:Exponential_revenue__Critical_w__Alpha_minus_1}. To get the explicit series expression in \refEquation{eqn:Exponential_revenue__Critical_w__Alpha_minus_1}, we integrate the generating function
\begin{equation}
\sum_{k=0}^\infty P_k(x) \varepsilon^k 
= ( 1 - 2x \varepsilon + \varepsilon^2 )^{ - \frac{1}{2} }
\end{equation}
of the Legendre polynomials over $x$ from $-1$ to $-1 - 2/t$, and we use for $k = 1, 2, \ldots$ that
\begin{equation}
P_{k+1}'(x) - P_{k-1}'(x) 
= (2k+1) P_k(x),
\quad
P_{k+1}(-1) - P_{k-1}(-1) = 0,
\end{equation}
see \cite[(4.7.29), (4.7.3-4)]{szego_orthogonal_1939} for the case $\lambda = 1/2$.

Case $\mathrm{(ii)}$. When $\alpha = 1/2$, we can write the threshold equation as
\begin{equation}
2(1-w)^{\frac{1}{2}} 
= 2 + \gamma B \varepsilon - (1+\gamma B) w.
\end{equation}
After squaring, we get the equation
\begin{equation}
w^2 + 2 \frac{ 2 - (2+\gamma B \varepsilon) ( 1 + \gamma B) }{ (1+\gamma B)^2 } w = \frac{ 4 - (2+\gamma B \varepsilon)^2 }{ (1+\gamma B)^2}.
\end{equation}
After a lengthy calculation, this yields the two solutions
\begin{equation}
w 
= \frac{2 \gamma B}{ (1+\gamma B)^2 } \Bigl( 1 + \frac{1}{2} ( 1 + \gamma B ) \varepsilon \pm \sqrt{ 1 - \frac{1+\gamma B}{\gamma B} \varepsilon } \Bigr).
\label{eqn:Exponential_revenue__Two_solutions_of_w__Alpha_minus_half}
\end{equation}
Noting that $-1 < \gamma B < 0$ in this case, and that $w$ is small positive when $\varepsilon$ is small positive, we take the $-$ sign in \refEquation{eqn:Exponential_revenue__Two_solutions_of_w__Alpha_minus_half}, and arrive at the square-root expression in \refEquation{eqn:Exponential_revenue__Critical_w__Alpha_half}, with $t$ given earlier. The series expansion given in \refEquation{eqn:Exponential_revenue__Critical_w__Alpha_half} and its validity range follow directly from this.

Case $\mathrm{(iii)}$. When $\alpha = 2$, we have $\gamma B > 0$, and the threshold equation can be written as
\begin{equation}
w^2 + 2 \gamma B w = 2 \gamma B \varepsilon.
\end{equation}
Using again that $w$ is small positive when $\varepsilon$ is small positive, the result in \refEquation{eqn:Exponential_revenue__Critical_w__Alpha_2} readily follows.
}

%%%%%%%%%%%%%%%%%%%%%%%%%%%%%%
%%% End of the appendices. %%%
%%%%%%%%%%%%%%%%%%%%%%%%%%%%%%
}

\myAcknowledgment{This research was financially supported by The Netherlands Organization for Scientific Research (NWO) in the framework of the TOP-GO program and by an ERC Starting Grant.}

\ifthenelse{\OperationsResearchTemplate = 1}
{\bibliographystyle{ormsv080}}
{\bibliographystyle{alpha}}

\bibliography{Bibliography}

\end{document}